\documentclass[a4paper,reqno,11pt]{amsart}
\usepackage{amssymb,amsfonts,verbatim,mathtools}
\usepackage[alwaysadjust]{paralist}
\usepackage[alwaysadjust]{paralist}
\numberwithin{equation}{section}
\newtheorem{notn}[equation]{Notation}
\newtheorem{claim}[equation]{Claim} %[section]
\newtheorem{cor}[equation]{Corollary}
\newtheorem{lem}[equation]{Lemma}
\newtheorem{prop}[equation]{Proposition}
\newtheorem{thm}[equation]{Theorem}
\theoremstyle{definition}
\newtheorem{defn}[equation]{Definition}

\newtheorem{rem}[equation]{Remark}

\newcommand{\cref}[1]{Corollary~\ref{#1}}

\begin{document}

%%%%%%%%%%%%%%%%%%%%%%%%%%%%%%%%%%%%%%%%%%%%%%

\title{On topological upper-bounds on the number of small cuspidal eigenvalues}
\author{Sugata Mondal}
\address{address 1}
%Max Planck Institute for Mathematics, Bonn} %, Germany}
%Vivatsgasse 7, 53111 Bonn, Deutschland}

\date{\today}

\subjclass{...}
\maketitle
\begin{abstract}
Let $S$ be a noncompact, finite area hyperbolic surface of type $(g, n)$. Let $\Delta_S$ denote the Laplace operator on $S$. As $S$ varies over the {\it moduli space} ${\mathcal{M}_{g, n}}$ of finite area hyperbolic surfaces of type $(g, n)$, we study, adapting methods of Lizhen Ji \cite{Ji} and Scott Wolpert \cite{Wo}, the behavior of {\it small cuspidal eigenpairs} of $\Delta_S$. In Theorem 2 we describe limiting behavior of these eigenpairs on surfaces ${S_m} \in {\mathcal{M}_{g, n}}$ when $({S_m})$ converges to a point in $\overline{\mathcal{M}_{g, n}}$. Then we consider the $i$-th {\it cuspidal eigenvalue}, ${\lambda^c_i}(S)$, of $S \in {\mathcal{M}_{g, n}}$. Since {\it non-cuspidal} eigenfunctions ({\it residual eigenfunctions} or {\it generalized eigenfunctions}) may converge to cuspidal eigenfunctions, it is not known if ${\lambda^c_i}(S)$ is a continuous function. However, applying Theorem 2 we prove that, for all $k \geq 2g-2$, the sets $${{\mathcal{C}_{g, n}^{\frac{1}{4}}}}(k)= \{ S \in {\mathcal{M}_{g, n}}: {\lambda_k^c}(S) > \frac{1}{4} \}$$ are open and contain a neighborhood of ${\cup_{i=1}^n}{\mathcal{M}_{0, 3}} \cup {\mathcal{M}_{g-1, 2}}$ in $\overline{\mathcal{M}_{g, n}}$. Moreover, using topological properties of nodal sets of {\it small eigenfunctions} from \cite{O}, we show that ${{\mathcal{C}_{g, n}^{\frac{1}{4}}}}(2g-1)$ contains a neighborhood of ${\mathcal{M}_{0, n+1}} \cup {\mathcal{M}_{g, 1}}$ in $\overline{\mathcal{M}_{g, n}}$. These results provide evidence in support of a conjecture of Otal-Rosas \cite{O-R}.
\end{abstract}
\section{Introduction}
In this paper a {\it hyperbolic surface} is a two dimensional complete Riemannian manifold $S$ with sectional curvature equal to $-1$. Such a surface is isomorphic to the quotient $\mathbb H/ \Gamma $, of the Poincar\'{e} upper halfplane $\mathbb H$ by a {\it Fuchsian group} $\Gamma$, i.e. a discrete torsion-free subgroup of PSL(2,$\mathbb R$). The {\it Laplace operator} on $\mathbb H$ is the differential operator which associates to a $C^2$-function $f$ the function $$\Delta f(z) = {y^2}(\frac{{\partial^2}f}{{\partial x}^2} + \frac{{\partial^2}f}{{\partial y}^2}).$$ Since the action of PSL(2, $\mathbb{R}$) on $\mathbb{H}$ leaves $\Delta$ invariant, $\Delta$ induces a differential operator on $S = {\mathbb H/ \Gamma}$ which extends to a self-adjoint operator $\Delta_S$ densely defined on ${L^2}(S).$ It is a general fact that the Laplace operator is a non-positive operator whose spectrum is contained in the smallest interval $( -\infty, -{\lambda_0}(S) ] \subset {\mathbb R^{-}} \cup \{0\}$ with ${\lambda_0}(S) \geq 0$.
\begin{defn}\label{d1}
Let $\lambda > 0$ be a real number and $f \in {L^2}(S)$ be a nonzero function on $S$. The pair $(\lambda, f)$ is called an eigenpair of $S$ if $ {\Delta_S} f + \lambda f \equiv 0 $ on $S$ where $\lambda$ and $f$ are respectively called an {\it eigenvalue} and an {\it eigenfunction} (sometimes a {\it $\lambda$-eigenfunction}). When $0 < \lambda \leq {1/4}$, we add the adjective {\it small} i.e. $(\lambda, f)$, $\lambda$ and $f$ are respectively called a {\it small eigenpair}, a {\it small eigenvalue} and a {\it small eigenfunction.}
\end{defn}
We begin with a noncompact, finite area hyperbolic surface $S$ of type $(g, n)$ i.e. $S \in {\mathcal{M}_{g, n}}$. The Laplace spectrum of such a surface is composed of two parts: {\it the discrete part} and {\it the continuous part} \cite{I}. The continuous part covers the interval $[\frac{1}{4}, \infty)$ and is spanned by {\it Eisenstein series} with multiplicity $n$. Eisenstein series are not eigenfunctions although they satisfy $$\Delta E(., s) + s(1-s) E(., s)=0 ,$$ because they are not in $L^2$. For this reason, they are called {\it generalized eigenfunctions}. The discrete spectrum consists of eigenvalues. They are distinguished into two parts: {\it the residual spectrum} and {\it the cuspidal spectrum}. An eigenpair $(\lambda, f)$ is called {\it residual} if $f$ is a linear combination of residues of meromorphic continuations of Eisenstein series. Such $\lambda$ and $f$ are respectively called a {\it residual eigenvalue} and a {\it residual eigenfuction}. The residual spectrum is a finite set contained in $[0, \frac{1}{4})$. On the other hand, an eigenpair $(\lambda, f)$ is called {\it cuspidal} if $f$ tends to zero at each cusp. In this case $\lambda$ and $f$ are respectively called a {\it cuspidal eigenvalue} and a {\it cuspidal eigenfuction}. These eigenvalues with multiplicity are arranged by increasing order and we denote ${\lambda^c_n}(S)$ the $n$-th cuspidal eigenvalue of $S$. For an arbitrary Fuchsian group $\Gamma$, it is not known whether the cardinality of the set of cuspidal eigenvalues of $\mathbb{H}/\Gamma$ is infinite. However, a famous theorem of A. Selberg says that it is the case when $\Gamma$ is arithmetic. Any cuspidal eigenpair $(\lambda, f)$ with $\lambda \leq \frac{1}{4}$ is called a {\it small cuspidal eigenpair} and in that case, $\lambda$ and $f$ are respectively called a {\it small cuspidal eigenvalue} and a {\it small cuspidal eigenfunction}.

In \cite{O-R}, Jean-Pierre Otal and Eulalio Rosas proved that the total number of small eigenvalues of any hyperbolic surface of type $(g, n)$ is at most $2g-3 +n$. In the same paper they formulate the following:

\textbf{Conjecture.} {\it Let $S$ be a noncompact, finite area hyperbolic surface of type $(g, n)$. Then ${\lambda^c_{2g-2}}(S) >  \frac{1}{4}$.}

This conjecture is motivated by the following two results
\begin{prop}\textbf{(Huxley \cite{Hu}, Otal \cite{O})}
Let $S$ be a finite area hyperbolic surface of genus $0$ or $1$. Then $S$ does not carry any small cuspidal eigenpair.
\end{prop}
\begin{prop}\textbf{(Otal \cite{O})}
Let $S$ be a finite area hyperbolic surface of type $(g, n)$. Then the multiplicity of a small cuspidal eigenvalue of $S$ is at most $2g-3$.
\end{prop}
The set ${\mathcal M}_{g,n}$ carries a topology for which two surfaces $\mathbb{H}/ {\Gamma}$ and $\mathbb{H}/ {\Gamma^{'}}$ are close when the groups $\Gamma$ and ${\Gamma^{'}}$ can be conjugated inside PSL$(2, \mathbb{R})$ so that they have generators which are close. With this topology ${\mathcal M}_{g,n}$ is not compact. However it can be compactified by adjoining ${\cup_i}{{\mathcal M}_{{g_i}, {n_i}}}$'s for each $({g_1}, {n_1}),..., ({g_k}, {n_k})$ with $2{\sum_i^k}({g_i} -2) + {\sum_i^k}{n_i} = 2g-2 +n$. In this compactification a sequence $(S_m) \in {{\mathcal M}_{g, n}}$ converges to ${S_\infty} \in \overline{{\mathcal M}_{g, n}}$ if and only if for any given $\epsilon > 0$ the $\epsilon$-thick part $({S_m^{[\epsilon, \infty)}})$ converges to ${S_\infty^{[\epsilon, \infty)}}$ in the Gromov-Hausdorff topology. Recall that the $\epsilon$-thick part of a surface $S$ is the subset of those points of $S$ where the {\it injectivity radius} is at least $\epsilon$. Recall also that the injectivity radius of a point $p \in S$ is the radius of the largest geodesic disc that can be embedded in $S$ with center $p$.

For any $N \in \mathbb{N}$ and $t \in \mathbb{R}_{>0}$ we define the sets $${\mathcal{C}^t_{g, n}}(N) = \displaystyle \{ S \in {\mathcal{M}_{g,n}}: {\lambda^c_N}(S) >  t \}.$$ It is clear that ${{\mathcal{C}_{g, n}^{\frac{1}{4}}}}(k) \subset {{\mathcal{C}_{g, n}^{\frac{1}{4}}}}(k+1)$ for $k \geq 1$. With this notation the conjecture can be formulated by saying that $${{\mathcal{C}_{g, n}^{\frac{1}{4}}}}(2g-2) = {\mathcal{M}_{g,n}}.$$ In this paper, we study the sets ${\mathcal{C}^\frac{1}{4}_{g, n}}(k)$. The methods developed here are not sufficient to prove the conjecture but we show that the sets ${\mathcal{C}^\frac{1}{4}_{g, n}}(2g-2)$ and ${\mathcal{C}^\frac{1}{4}_{g, n}}(2g-1)$ (${\mathcal{C}^\frac{1}{4}_{g, n}}(2g-2) \subseteq {\mathcal{C}^\frac{1}{4}_{g, n}}(2g-1)$) contains neighborhoods of certain strata in the compactification of $\mathcal{M}_{g, n}$.

\begin{thm}
$\textbf{(i)}$ For any integer $k$, ${{\mathcal{C}_{g, n}^{\frac{1}{4}}}}(k)$ is an open subset of ${\mathcal{M}_{g,n}}$.\\* \textbf{(ii)} ${{\mathcal{C}_{g, n}^{\frac{1}{4}}}}(2g-2)$ contains a neighborhood of ${\cup_{i=1}^n}{\mathcal{M}_{0, 3}} \cup {\mathcal{M}_{g-1, 2}}$ in $\overline{\mathcal{M}_{g, n}}$.\\* \textbf{(iii)} ${{\mathcal{C}_{g, n}^{\frac{1}{4}}}}(2g-1)$ contains a neighborhood of ${\mathcal{M}_{0, n+1}} \cup {\mathcal{M}_{g, 1}}$ in $\overline{\mathcal{M}_{g, n}}$.
\end{thm}
Observe that it is theoretically possible for a residual eigenfunction to converge to a cuspidal eigenfunction. Therefore indicating that ${\lambda^c_{2g-1}}$ may not be continuous. Also, the result \cite{P-S} suggest that ${\lambda^c_{2g-1}}$ may not be continuous at those $S \in {\mathcal{M}_{g, n}}$ where it takes value strictly more than $\frac{1}{4}$. Therefore, the first assertion is not completely trivial.

The paper is organized as follows. In \S1 we recall some preliminaries for convergence of hyperbolic surfaces in $\overline{\mathcal{M}_{g, n}}$. In \S2 and \S3 we study convergence properties of eigenpairs on converging hyperbolic surfaces. Similar study has already been carried out by Scott Wolpert \cite{Wo}, Lizhen Ji \cite{Ji} and Christopher Judge \cite{J}. We shall first make precise the notions of convergence in $\overline{{\mathcal M}_{g,n}}$ and the notion of convergence of a sequence of functions on a converging sequence of surfaces.
\subsection{Convergence of functions}
Let $({S_m})$ be a sequence of surfaces in ${\mathcal{M}_{g, n}}$ converging to a surface $S_\infty$ in the compactification $\overline{{\mathcal{M}_{g, n}}}$. Another way of understanding this convergence is as follows:

Let ${S_m} = {\mathbb{H}/ {\Gamma_m}}$ and let $0< {c_0} < {\epsilon_0}$ ($\epsilon_0$ is the Margulis constant; see thick/thin decomposition for details) be a fixed constant. Let ${x_m} \in {S^{[{c_0}, \infty)}_m}$. Up to a conjugation of $\Gamma_m$ in PSL$(2, \mathbb{R})$, one may assume that $i \in \mathbb{H}$ is mapped to $x_m$ under the projection $\mathbb{H} \rightarrow {\mathbb{H}/ {\Gamma_m}}$. Then up to extracting a subsequence we may suppose that $\Gamma_m$ converges to some Funchsian group $\Gamma_\infty$. We say that the pair $({\mathbb{H}/ {\Gamma_m}}, {x_m})$ converges to $({\mathbb{H}/ {\Gamma_\infty}}, {x_\infty})$ where ${x_\infty}$ is the image of $i \in \mathbb{H}$ under the projection $\mathbb{H} \rightarrow {\mathbb{H}/ {\Gamma_\infty}}$. Let $S_\infty$ be the hyperbolic surface of finite area whose connected components are the ${\mathbb{H}/ {\Gamma_\infty}}$'s for different choices of base point ${x_m}$ in different connected components of ${S^{[{c_0}, \infty)}_m}$. The surface ${S_\infty}$ does not depend, up to isometry, on the choice of the base point $x_m$ in a fixed connected component of ${S^{[{c_0}, \infty)}_m}$ (i.e. if ${y_m}$ be a point in the same connected component of ${S^{[{c_0}, \infty)}_m}$ as $x_m$ then the corresponding limiting surfaces are isometric). One can check that $({S_m}) \rightarrow {S_\infty}$ in $\overline{{\mathcal{M}_{g, n}}}$.

\underline{\textbf{Convergence of functions}}\\*
Fix an $\epsilon > 0$ and choose a base point ${x_m} \in {{S_m}^{[\epsilon, \infty)}}$ for each $m$. Assume that the pair $({\mathbb{H}/ {\Gamma_m}}, {x_m})$ converges to $({\mathbb{H}/ {\Gamma_\infty}}, {x_\infty})$ where, for each $m \in \mathbb{N} \cup \{ \infty \}$, the point $i \in \mathbb{H}$ maps to $x_m$ under the projection $\mathbb{H} \rightarrow {\mathbb{H} / {\Gamma_m}}$.

For a $C^\infty$ function $f$ on $S_m$ denote by $\widetilde{f}$ the lift of $f$ under the projection $\mathbb{H} \rightarrow {\mathbb{H}/ {\Gamma_m}}$. Let $({f_m})$ be a sequence of functions in ${C^\infty}({S_m}) \cap {L^2}({S_m})$. One says that $({f_m})$ {\it converges to} a continuous function $f_\infty$ if $\widetilde{{f_m}}$ converges, uniformly over compact subsets of $\mathbb{H}$, to $\widetilde{{f_\infty}}$ for each choice of base points ${x_m} \in {{S_m}^{[\epsilon, \infty)}}$ and for each $\epsilon < {\epsilon_0}$.

With the above understanding of convergence of functions we shall prove the following theorem which has close resemblance with \cite[Theorem 1.2]{Ji} and \cite[Theorem 4.2]{Wo}. However, our result does not follow from these. We would like to mention that a similar limiting theorem might not be true (see \cite[p-71]{Wo} if one considers ${\lambda_m} \geq \frac{1}{4}$ instead of ${\lambda_m} \leq \frac{1}{4}$ (see Theorem 2).

In the following, for a function $f \in {L^2}(S)$, we shall denote the ${L^2}$ norm of $f$ by $\Vert f \rVert$. Also, for $f \in {L^2}(V)$ and $U \subset V$ we denote the $L^2$-norm of the restriction of $f$ to $U$ by ${{\lVert f \rVert}_U}$. A function $f \in {L^2}(V)$ will be called {\it normalized} if $\lVert f \rVert =1$. An eigenpair $(\lambda, \phi)$ will be called normalized if $\phi$ is normalized.
\begin{thm}
Let ${S_m}\rightarrow {S_\infty}$ in $\overline{{\mathcal{M}_{g, n}}}$. Let $({\lambda_m}, {\phi_m})$ be a normalized small cuspidal eigenpair of $S_m$. Assume that ${\lambda_m}$ converges to ${\lambda_\infty}$. Then one of the following holds:\\*
$(1)$ There exist strictly positive constants $\epsilon, \delta$ such that $\limsup{{\lVert {\phi_m} \rVert}_{S_m^{[\epsilon, \infty)}}} \geq \delta$. Then, up to extracting a subsequence, $({\phi_m})$ converges to a $\lambda_\infty$-eigenfunction $\phi_\infty$ of $S_\infty$.\\*
$(2)$ For each $\epsilon >0$ the sequence $({{\lVert {\phi_m} \rVert}_{S_m^{[\epsilon, \infty)}}}) \rightarrow 0.$ Then ${S_\infty} \in \partial {\mathcal{M}_{g, n}}$ and ${\lambda_\infty}= \frac{1}{4}$. Moreover, there exist constants ${K_m} \rightarrow \infty$ such that, up to extracting a subsequence, $({K_m}{\phi_m})$ converges to a linear combination of Eisenstein series and (possibly) a cuspidal ${\lambda_\infty}$-eigenfunction of $S_\infty$.
\end{thm}
\begin{rem}
For $s=\frac{1}{2}$, by Eisenstein series we understand a linear combination of the following two:\\* \textbf{(i)} the classical (meromorphic continuation) Eisenstein series ${E^i}(., \frac{1}{2})$ corresponding to the cusps ($i$ is the index for cusps) on the surface,\\* \textbf{(ii)} the derivatives ${\frac{\partial}{\partial s}{E^i}(., s)}|_{s= \frac{1}{2}}$ of ${E^i}(., s)$ at $s = \frac{1}{2}$.\\* The first Fourier coefficient of such functions in any cusp have the form $\alpha {y^\frac{1}{2}} + \beta {y^\frac{1}{2}}\log{y}$. Each moderate growth $\frac{1}{4}$-eigenfunction is a linear combination of Eisenstein series, in the above sense, and (possibly) a cuspidal eigenfunction.
\end{rem}
Theorem 2 will be applied to prove all three statements of Theorem 1. The first one is a direct application; in \S4 we prove:

\underline{\textbf{Lemma 1}} {\it For any $k \geq 1$, ${\mathcal{C}_{g, n}^{\frac{1}{4}}}(k)$ is an open subset of $\mathcal{M}_{g, n}$.}

The second statement of Theorem 1 is also an easy application of Theorem 2 and the Buser construction \cite{Bu}: we explain it now since the proof is short. We argue by contradiction and assume that there is a sequence $({S_m})$ in $\mathcal{M}_{g, n}$ such that $S_m$ converges to ${S_\infty} \in {\cup_{i=1}^n}{\mathcal{M}_{0, 3}} \cup {\mathcal{M}_{g-1, 2}}$ and ${\lambda_{2g-2}^c}({S_m}) \leq \frac{1}{4}$. Then $S_\infty$ has exactly $n + 1$ components of which exactly $n$ are thrice punctured spheres. Observe that each component of $S_\infty$ contains an {\it old cusp} i.e. cusps of $S_\infty$ which are limits of cusps of $S_m$ (see Proof of Theorem 2).

The construction used in the proof of \cite[Theorem 8.1.3]{Bu} implies that, for $m$ large, $S_m$ has at least $n$ eigenvalues that converge to zero as $m$ tends to infinity. Let us suppose by contradiction that one of the corresponding eigenfunctions $\phi_m$ is cuspidal. Then by Theorem 2, $\phi_m$ converges uniformly over compacta to a function $\phi$ and $\phi$ is an eigenfunction for the eigenvalue $0$. So $\phi$ is constant in each component of $S_\infty$. On those components of $S_\infty^{[\epsilon, \infty)}$ that contains an old cusp $\phi$ is necessarily zero because $\phi_m$ being cuspidal the average of $\phi_m$ over any horocycle is zero. On the the other component (the one that does not contain an old cusp) $\phi$ is zero because the mean of $\phi$ over $S_\infty$ is equal to the mean of $\phi_m$ over $S_m$ which is zero (follows from Theorem \ref{H}). Therefore, $\phi$ is the zero function which is a contradiction by Theorem 2. Hence, for large $m$ each eigenfunction corresponding to any of the first $n$ eigenvalues of $S_m$ is necessarily residual. Now if ${\lambda_{2g-2}^c}({S_m}) \leq \frac{1}{4}$ then each $S_m$ has at least $2g-2+n$ small eigenvalues. This is a contradiction to \cite[Theorem 2]{O-R}. Therefore we have proved that ${\mathcal{C}_{g, n}^{\frac{1}{4}}}(2g-2)$ contains a neighborhood of ${\cup_{i=1}^n}{\mathcal{M}_{0, 3}} \cup {\mathcal{M}_{g-1, 2}}$ in $\overline{\mathcal{M}_{g, n}}$.

In the last section we prove the last statement of Theorem 1. We consider ${\mathcal{M}_{g, 1}} \cup {\mathcal{M}_{0, n+1}}$ as a subset of $\partial{\mathcal{M}_{g, n}}= \overline{\mathcal{M}_{g, n}} \setminus {\mathcal{M}_{g, n}}$ and show the following

\begin{prop}\label{neighbrhd}
There exists a neighborhood $\mathcal{N}({\mathcal{M}_{g, 1}} \cup {\mathcal{M}_{0, n+1}})$ of ${\mathcal{M}_{g, 1}} \cup {\mathcal{M}_{0, n+1}}$ in $\overline{\mathcal{M}_{g, n}}$ such that for each $S \in \mathcal{N}({\mathcal{M}_{g, 1}} \cup {\mathcal{M}_{0, n+1}})$: ${\lambda^c_{2g-1}}(S) > \frac{1}{4}$ i.e. $$\mathcal{N}({\mathcal{M}_{g, 1}} \cup {\mathcal{M}_{0, n+1}}) \subset {{\mathcal{C}_{g, n}^{\frac{1}{4}}}}(2g-1).$$
\end{prop}
Now we briefly sketch a proof of this proposition. We argue by contradiction and consider a sequence $({S_m})$ in $\mathcal{M}_{g, n}$ that converges to $S_\infty$ in $({\mathcal{M}_{g, 1}} \cup {\mathcal{M}_{0, n+1}}) \subset \partial{\mathcal{M}_{g, n}}$ such that ${{\lambda^c}_{2g-1}}({S_m}) \leq  \frac{1}{4}$. Then, for $1 \leq i \leq 2g-1$ and for each $m$, we choose a small cuspidal eigenpair $({\lambda^i_m}, {\phi^i_m})$ of $S_m$ such that $\\(i)~ \{{\phi^i_m}\}_{i=1}^{2g-1}$ is an orthonormal family in ${L^2}({S_m}),$ $\\(ii)~ {\lambda_m^i}$ is the $i$-th eigenvalue of $S_m$.

For $1 \leq i \leq 2g-1$ let $({\lambda^i_m})$ converges to ${\lambda^i_\infty}$ as $m \rightarrow \infty$. By Theorem 2 there are two possible types of behavior that the sequence $({\phi^i_m})$ can exhibit. Either, for each $1 \leq i \leq 2g-1$ the sequence $({\phi^i_m})$ converges to a ${\lambda^i_\infty}$-eigenfunction ${\phi^i_\infty}$ on $S_\infty$, or for some $i$ the sequence $({\lambda^i_m}, {\phi^i_m})$ satisfies condition $(2)$ in Theorem 2. However, in our case we have the following lemma:

\underline{\textbf{Lemma 2}}~ {\it For each $i$, $1 \leq i \leq 2g-1$, up to extracting a subsequence, the sequence $({\phi^i_m})$ converges to a $\lambda^i_\infty$-eigenfunction $\phi^i_\infty$ of $S_\infty$. The limit functions ${\phi^i_\infty}$ and ${\phi^j_\infty}$ are orthogonal for $i \neq j$ i.e. $S_\infty$ has at least $2g-1$ small eigenvalues. Moreover none of the ${\phi^i_\infty}$ is residual.}

Then we count the number of small eigenvalues of $S_\infty$ using \cite{O-R} to conclude that at least one of ${\phi^i_\infty}$ is nonzero on the component of $S_\infty$ of type $(0, n+1)$. This leads to a contradiction by Huxley \cite{Hu} or \cite[Proposition 2]{O}.
\subsection{Acknowledgement}
 The author would like to express his sincere gratitude to his advisor Jean-Pierre Otal for his patience, encouragement and insight. The author was supported during this research by the Indo-French CEFIPRA-IFCPAR grant.
\section{Preliminaries}
In this section we shall recall some preliminary concepts that are important for our purpose. Metric convergence of a sequence $(S_m) \in {{\mathcal M}_{g, n}}$ to ${S_\infty} \in \overline{{{\mathcal M}_{g, n}}}$ is one of the prime aspects of our study. We start by explaining the thick/thin decomposition of a hyperbolic surface which is convenient to understand the metric convergence.
\subsection{The thick / thin decomposition of a hyperbolic surface}
Let $S \in {\mathcal{M}_{g, n}}.$ Recall that for any $\epsilon >0$, {\it the $\epsilon$-thin part of $S$}, $S^{(0, \epsilon)}$, is the set of points of $S$ with injectivity radius $< \epsilon$. The complement of $S^{(0, \epsilon)}$, {\it the $\epsilon$-thick part of $S$}, denoted by $S^{[\epsilon, \infty)}$, is the set of points where the injectivity radius of $S$ is $\geq \epsilon$.
\subsubsection{Cylinders}\label{cylinder}
Let $\gamma$ be a simple closed geodesic on $S$. It can be viewed as the quotient of a geodesic in $\mathbb{H}$ by a hyperbolic isometry $\Upsilon$ fixing the geodesic. We may conjugate $\Upsilon$ such that the geodesic is the imaginary axis and the isometry is $\tau: z \rightarrow {e^{2\pi l}}z$, $2\pi l= {l_\gamma}$ being the length of the geodesic. We define the {\it hyperbolic cylinder $~\mathcal{C}$ with core geodesic $\gamma$} as the quotient ${\mathbb{H}/ <\tau>}$. Recall that the {\it Fermi coordinates} on $\mathcal{C}$ assign to each point $p \in \mathcal{C}$ the pair $(r, \theta) \in {\mathbb{R}} \times \{ \gamma \}$ where $r$ is the signed distance of $p$ from $\gamma$ and $\theta$ is the projection of $p$ on $\gamma$ \cite[p. 4]{Bu}. These coordinates give a diffeomorphism of this hyperbolic cylinder to $\mathbb{R} \times {\mathbb{R}/{2\pi}\mathbb{Z}}$. In terms of these coordinates the hyperbolic metric is given by: $$d{s^2} = d{r^2} + {l^2} {\cosh^2}{r} d{\theta^2}.$$ For $w \geq l$ we define the {\it collar} $~{\mathcal{C}^w}$ around $\gamma$ by $${\mathcal{C}^w}= \{(r, \theta) \in \mathcal{C}: {l_\gamma}\cosh r < w , 0 \leq \theta \leq 2\pi \}.$$ Then ${\mathcal{C}^w}$ is diffeomorphic to an annulus whose each boundary component has length $w$. The {\it Collar Theorem} of Linda Keen \cite{Ke} says that $\mathcal{C}^1$ embeds in $S$ (more precisely, $\mathcal{C}^{w({l_\gamma})}$ embeds in $S$ where $w({l_\gamma}) = {l_\gamma} \cosh ({\sinh^{-1}}(\frac{1}{\sinh{\frac{l_\gamma}{2}}})) > 1$ and $w({l_\gamma}) \approx 2$).
\subsubsection{Cusps}\label{cuspd}
$S$ has $n$ ends called {\it punctures}. {\it Cusps} are particular neighborhood of the punctures. Denote by $\iota$ the parabolic isometry $\iota: z \rightarrow z + 2\pi$. For a choice of $t > 0$, a cusp $\mathcal{P}^t$ is the half-infinite cylinder $\{z= x+iy : y > \frac{2\pi}{t} \}/<\iota>$. The boundary curve $\{y = \frac{2\pi}{t} \}$ is a {\it horocycle} of length $t$ that we identify with $\mathbb{R}/ t\mathbb{Z}$. One can parametrize $\mathcal{P}^t$ using the {\it horocycle coordinates} \cite[p. 4]{Bu} with respect to its boundary horocycle $\{y = \frac{2\pi}{t} \}$. The horocycle coordinates assigns to a point $p \in {\mathcal{P}^t}$ the pair $(r, \theta) \in {\mathbb{R}_{\geq 0}} \times \{ \mathbb{R}/ t\mathbb{Z} \}$ where $r$ is the distance from $p$ to the horocycle and $\theta$ the projection of $p$ on the horocycle. In terms of these coordinates the hyperbolic metric takes the form: $$d{s^2} = d{r^2} + {(\frac{t}{2\pi})^2}{e^{-2r}} d{\theta^2}.$$

Recall that the cusp $\mathcal{P}^1$ (in fact $\mathcal{P}^2$) around each puncture embeds in $S$ and that those cusps corresponding to distinct punctures have disjoint interiors (ref. \cite[Chapter 4]{Bu}). We call them {\it standard cusp}s. Observe that the area and boundary length of a standard cusp is equal to $1$. For $t\leq 1$ denote the disjoint union ${\bigcup_{c \in S}} ~{\mathcal{P}^t}$ by ${S_c^{(0, t)}}$ where $c$ ranges over distinct cusps in $S$.
\subsubsection{The decomposition}
By Margulis lemma there exists a constant ${\epsilon_0} >0$, the Margulis constant, such that for all $\epsilon \leq {\epsilon_0}$, $S^{(0,\epsilon)}$ is a disjoint union of embedded collars, one for each geodesic of length less than $2\epsilon$, and of embedded cusps, one for each puncture. The collar around a geodesic of length $\leq \epsilon$ is called a {\it Margulis tube}.
\subsection{Metric degeneration of a collar to a pair of cusps}
We describe how a collar around a geodesic of length ${l_\gamma}= 2\pi l$ converges as $l$ tends to zero to a pair of cusps. First shift the origin of the Fermi coordinates of ${\mathcal{C}^{w({l_\gamma})}}$ to the right boundary of ${\mathcal{C}^{w({l_\gamma})}}$ by making the change of variable $t= r - {\sinh^{-1}}(\frac{1}{\sinh{\frac{l_\gamma}{2}}})$. In the shifted Fermi coordinates the metric on ${\mathcal{C}^{w({l_\gamma})}}$ is equal to $$d{s^2} = d{r^2} + {l^2}{\cosh^2}(r - {\sinh^{-1}}(\frac{1}{\sinh{\frac{l_\gamma}{2}}}))d{\theta^2}.$$ For $r$ in a compact region we have the limiting $${\lim_{l \rightarrow 0}} ~{l {\cosh(r- {\sinh^{-1}}(\frac{1}{\sinh{\frac{l_\gamma}{2}}}))}} = \frac{e^{-r}}{\pi}.$$ Now the hyperbolic metric on $\mathcal{P}^2$ is equal to $${d{s^2}}= d{r^2} + \frac{e^{-2r}}{\pi^2}d{\theta^2}$$ with respect to the boundary horocycle $\{ y = \pi \}$ of $\mathcal{P}^2$.

Choose a base point $p_l$ on the right half of ${\mathcal{C}^{w({l_\gamma})}}^{[\epsilon, \infty)}$. Then by above, as $l \rightarrow 0$, the pair $({\mathcal{C}^{w({l_\gamma})}}, {p_l})$ converges, up to extracting a subsequence, to $({\mathcal{P}^2}, p)$ where $p \in {\mathcal{P}^2}^{[\epsilon, \infty)}$. Since one can choose the base point on the left half of ${\mathcal{C}^{w({l_\gamma})}}$ also, the metric limit of ${\mathcal{C}^{w({l_\gamma})}}$ is a pair of $\mathcal{P}^2$.
\section{Mass distribution of small cuspidal functions over thin parts}
Our goal is to study the behavior of sequences of small cuspidal eigenpairs $({\lambda_n}, {f_n})$ of ${S_n} \in {\mathcal{M}_{g, n}}$ when $(S_n)$ converges to ${S_\infty} \in \overline{{{\mathcal M}_{g, n}}}$ and finally to prove Theorem 1. For this we need to understand how the {\it mass} (${L^2}$ norm) of a small eigenfunction is distributed over the surface, and in particular how it is distributed with respect to the {\it thin/thick} decomposition. Let $S \in {\mathcal{M}_{g, n}}$. Recall that for any $\epsilon \leq {\epsilon_0}$ the $\epsilon$-thin part, $S^{(0, \epsilon)}$, of $S$ consists of cusps and Margulis tubes. We separately study the mass distribution of a small cuspidal eigenfunction over these two different types of domains.
\subsection{Mass distribution over cusps}
For $2\pi \leq a < b$ consider the annulus ${\mathcal{P}(a, b)} = \{ (x,y) \in {\mathcal{P}^1}: a \leq y < b \}$ contained in a cusp $\mathcal{P}^1$ and bounded by two horocycles of length $\frac{2\pi}{a}$ and $\frac{2\pi}{b}$. We begin our study with the following lemma.
\begin{lem}\label{cusp}
For any $b>2\pi$ there exists $K(b)< \infty$ such that for any small cuspidal eigenpair $(\lambda, f)$ of $\mathcal{P}^1$ one has
\begin{equation}
{{\lVert f \rVert}_{\mathcal{P}(b,\infty)}} < K(b) {{\lVert f \rVert}_{\mathcal{P}(2\pi, b)}}.
\end{equation}
If $\lambda < \frac{1}{4} - \eta$ for some $\eta > 0$ then there exists a constant $T(b, \eta) < \infty$ depending on $b$ and $\eta$ such that for any small eigenpair $(\lambda, f)$ one has
\begin{equation}
{{\lVert f \rVert}_{\mathcal{P}(b,\infty)}} < T(b, \eta) {{\lVert f \rVert}_{\mathcal{P}(2\pi, b)}}.
\end{equation}
Furthermore, $K(b), T(b, \eta) \rightarrow 0$ as $b \rightarrow \infty$.
\end{lem}
\textbf{Proof.} We begin with the first part. Since $f$ is cuspidal inside $\mathcal{P}^1$ it can be expressed as
\begin{equation}\label{expres}
 f(z) = \sum_{n \in {{\mathbb Z}^*}} {f_n} {W_s}(nz)
\end{equation}
where $s(1-s) = \lambda$ and $W_s$ is the {\it Whittaker function} (see \cite[Proposition 1.5]{I}). The meaning of (\ref{expres}) is that the right hand series converges to $f$ in ${L^2}({\mathcal{P}^1})$ and that the convergence is uniform over compact subsets. Recall also that for $n \in {\mathbb Z}^*$ the Whittaker functions is defined by $${W_s}(nz) =  2{(|n|y)^{\frac{1}{2}}} {K_{s-\frac{1}{2}}}(|n| y) {e^{inx}}$$ where $K_\epsilon$ is the {\it McDonald's function} and that for any $\epsilon$ (see \cite[p. 119]{Le})
\begin{equation}
{{K_\epsilon}(y)}= \displaystyle \frac{1}{2} {\int_{-\infty}^{+\infty} e^{ -y \cosh {u} - \epsilon u} du}
\end{equation}
whenever the integral makes sense. From the expression it is clear that the functions $({W_s}(n . ))$ form an orthogonal family over $\mathcal{P}(a, b)$ (independent of the choices of $a$ and $b$). Hence (1) will follow from the following claim.
\begin{claim}\label{cusp3}
Let $s \in [\frac{1}{2}, 1]$. Then for any $b > 2\pi$ there exists $K(b) < \infty$ such that for all $n \in {{\mathbb Z}^*}$ $${{\lVert {W_s}(nz) \rVert}_{\mathcal{P}(b,\infty)}} \leq K(b) {{\lVert {W_s}(nz) \rVert}_{\mathcal{P}(2\pi, b)}}.$$ Furthermore, $K(b) \rightarrow 0$ as $b \rightarrow \infty$.
\end{claim}
\textbf{Proof.} From the expression of $W_s$ we have $${{\lVert {{W_s}(nz)} \rVert}_{\mathcal{P}(a, b)}} = 2\pi \bigg({\int_a^b} 4|n|y {{K_{s-\frac{1}{2}}}(|n| y)}^2 \frac{dy}{y^2}\bigg).$$ To prove the claim we may suppose that $n \geq 1$. Our next objective is to obtain bounds for the functions ${{K_{s-\frac{1}{2}}}(y)}$ for $s \in [\frac{1}{2}, 1]$. We start from the above integral representation of ${{K_\epsilon}(y)}$. We write ${{K_\epsilon}(y)} = \displaystyle \frac{1}{2} \{ c(\epsilon, y) + d(\epsilon, y) \}$ where
\begin{equation}
 c(\epsilon, y) = \displaystyle \int_{-1}^{1} e^{-y\cosh u - \epsilon u} du
\end{equation}
and
\begin{equation}
 d(\epsilon, y) =\displaystyle \int_{-\infty}^{-1} e^{-y\cosh u - \epsilon u} du + \int_1^{\infty} e^{-y\cosh u - \epsilon u} du.
\end{equation}
Now we treat $c(\epsilon, y)$ and $d(\epsilon, y)$ separately.

\underbar{Bounding $c(\epsilon, y)$}:

We have $$c(\epsilon, y)= \int_{-1}^1 e^{-y\cosh u}. e^{-\epsilon u} du \leq e^{\epsilon}. \int_{-1}^1 e^{-y \cosh u} du = e^{\epsilon}\int_{-1}^1 e^{-y(1+ \frac{u^2}{2!} + \frac{u^4}{4!} +...)} du$$ $$= e^{\epsilon} . e^{-y} \int_{-1}^1 e^{-y( \frac{u^2}{2!} + \frac{u^4}{4!} + ...)} du \leq 2 e^{\epsilon} .e^{-y} \int_{0}^1 e^{-y \frac{u^2}{2!}} du.$$ Since $ e^{\frac{y u^2 }{2}} > 1 +\frac{ y u^2}{2} ~ \textrm{for} ~ u > 0 $, we have: $$ \int_{0}^1 e^{-y \frac{u^2}{2!}} du < \int_{0}^1 \displaystyle \frac{du}{1 + \frac{y u^2}{2}} = \frac{2}{y} {{\tan}^{-1}} (\frac{y}{2}) \leq \frac{2}{y}. \frac{\pi}{2}.$$ Therefore $$c(\epsilon, y) \leq 2 \pi e^{\epsilon} \frac{e^{-y}}{y}.$$ To obtain a lower bound, we write $$\int_{-1}^1 e^{-y\cosh u}. e^{-\epsilon u} du \geq e^{-\epsilon}. \int_{-1}^1 e^{-y \cosh u} du= e^{-\epsilon}\int_{-1}^1 e^{-y(1+ \frac{u^2}{2!} + \frac{u^4}{4!} +...)} du$$ $$= 2 e^{-\epsilon}. e^{-y} \int_{0}^1 e^{-y( \frac{u^2}{2!} + \frac{u^4}{4!} + ...)} du.$$ Since for all $u \in (0, 1]$ one has $$\displaystyle \frac{u^2}{2!} + \frac{u^4}{4!} + ... < u (\frac{1}{2} + \frac{1}{4} + \frac{1}{8} + ...) = u.$$ Hence $$c(\epsilon, y) \geq 2 e^{-\epsilon} . e^{-y} \int_0^1 e^{-uy} du = 2 e^{-\epsilon} \frac{e^{-y}}{y} (1 - e^{-y}).$$ Combining the above two inequalities $$2 e^{-\epsilon} \frac{e^{-y}}{y} (1 - e^{-y}) \leq c(\epsilon, y) \leq 2 \pi e^{\epsilon} \frac{e^{-y}}{y}.$$

\underbar{Bounding $d(\epsilon, y)$}:

$$d(\epsilon, y)= \int_{-\infty}^{-1} e^{-y\cosh u - \epsilon u} du + \int_1^{\infty} e^{-y\cosh u - \epsilon u} du $$ $$= \int_{1}^{\infty} e^{-y\cosh u - \epsilon u} du + \int_{1}^{\infty} e^{-y\cosh u + \epsilon u} du.$$ Now for any $u> 1$, $$\displaystyle \frac{u^2}{2!} + \frac{u^4}{4!} + ...  > {\gamma_0} u^2 > {\gamma_0} u$$ where ${\gamma_0} = \displaystyle {\sum_{n=1}^\infty} \frac{1}{(2n)!}$.

Thus $$d(\epsilon, y) = e^{-y} \int_1^{\infty} \{ e^{-y(\frac{u^2}{2!} + \frac{u^4}{4!} + ...) -\epsilon u} + e^{-y(\frac{u^2}{2!} + \frac{u^4}{4!} + ...) + \epsilon u} \} du  $$ $$\leq e^{-y} \int_1^{\infty} \{ e^{-y{\gamma_0} u - \epsilon u} + e^{-y{\gamma_0} u + \epsilon u} \} du $$ $$= \displaystyle \frac{e^{-y}}{y} \bigg( \frac{e^{-(y{\gamma_0} + \epsilon)}}{{\gamma_0} + \frac{\epsilon}{y}} + \frac{e^{-(y{\gamma_0} - \epsilon)}}{{\gamma_0} - \frac{\epsilon}{y} } \bigg).$$

Thus combining the estimates for $c(\epsilon, y)$ and $d(\epsilon, y)$ we obtain $$2 e^{-\epsilon} \frac{e^{-y}}{y} (1 - e^{-y}) < {K_{\epsilon}}(y) < 2 \pi e^{\epsilon} \frac{e^{-y}}{y} + \displaystyle \frac{e^{-y}}{y} \bigg( \frac{e^{-(y{\gamma_0} + \epsilon)}}{{\gamma_0} + \frac{\epsilon}{y}} + \frac{e^{-(y{\gamma_0} - \epsilon)}}{{\gamma_0} - \frac{\epsilon}{y} } \bigg).$$ Let $$\delta(\epsilon, y) =  \frac{e^{-(y{\gamma_0} + \epsilon)}}{{\gamma_0} + \frac{\epsilon}{y}} + \frac{e^{-(y{\gamma_0} - \epsilon)}}{{\gamma_0} - \frac{\epsilon}{y} }.$$ Observe that for $\epsilon < 1$ and $ y \geq \frac{2}{\gamma_0}$ $$\delta(\epsilon, y) < \frac{4 \cosh 1}{\gamma_0} {e^{-{\gamma_0}y}} = {\delta_0}(y).$$ So, for $ y \geq \frac{2}{\gamma_0}$ large enough
\begin{equation}
 2e^{-\epsilon} \frac{e^{-y}}{y} < {K_\epsilon}(y) < \frac{e^{-y}}{y} \bigg(2\pi e^{\epsilon} + {\delta_0}(y) \bigg).
\end{equation}
Going back to the expression of $W_s$, for $s \in [\frac{1}{2}, 1]$, we find: $$\frac{1}{2\pi}{{\lVert {W_s}(nz) \rVert}_{\mathcal{P}(2\pi, b)}^2} =  {\int_{2\pi}^b} 4ny {{K_{s-\frac{1}{2}}}(n y)}^2 \frac{dy}{y^2} =  {\int_{2\pi}^b} 4n {{{K_{s-\frac{1}{2}}}(n y)}^2} \frac{dy}{y}$$ $$\geq  {\int_{2\pi}^b} \displaystyle \frac{4n}{b} {{K_{s-\frac{1}{2}}}(n y)}^2 dy > \displaystyle \frac{16n {e^{1-2s}}}{b} {\int_{2\pi}^b} \frac{e^{-2ny}}{(ny)^2} dy = \displaystyle \frac{16n{e^{1-2s}}}{{n^2} b} {\int_{2\pi}^b} \frac{e^{-2ny}}{y^2} dy $$ $$=\displaystyle \frac{16n {e^{1-2s}}}{{n^2} b} \bigg( {\int_{2\pi}^\frac{b}{2}} \frac{e^{-2ny}}{y^2} dy + {\int_\frac{b}{2}^b} \frac{e^{-2ny}}{y^2} dy \bigg) > \displaystyle \frac{16n {e^{1-2s}}}{{n^2} b} \bigg({\int_\frac{b}{2}^b} \frac{e^{-2ny}}{y^2} dy \bigg)$$ $$= \frac{16 {e^{1-2s}}}{nb} \frac{e^{-nb}}{n \frac{b^2}{4}} \{ 1+ O( e^{-nb} + {\frac{2}{b}}) \}$$ i.e.
\begin{equation}\label{lb}
 {{\lVert {W_s}(nz) \rVert}_{\mathcal{P}(2\pi, b)}^2} > 2\pi \frac{16 {e^{1-2s}}}{nb} \frac{e^{-nb}}{n \frac{b^2}{4}} \{ 1+ O( e^{-nb} + {\frac{1}{b}}) \}
\end{equation}
Also, $$\frac{1}{2\pi} {{\lVert {W_s}(nz) \rVert}_{\mathcal{P}(b, \infty)}^2} = {\int_b^\infty} \displaystyle 4n y {{K_{s-\frac{1}{2}}}(n y)}^2 \frac{dy}{y^2} = {\int_b^{\infty}} 4n {{K_{s-\frac{1}{2}}}(n y)}^2 \frac{dy}{y}$$ $$\leq {\int_b^{\infty}} \displaystyle \frac{4n}{b} {{K_{s-\frac{1}{2}}}(n y)}^2 dy \leq \frac{4n(2\pi e^{(s-\frac{1}{2})} + {\delta_0}(b))^2}{b}{\int_b^{\infty}} \frac{e^{-2ny}}{(ny)^2} dy$$ $$= \frac{4{{(2\pi e^{(s-\frac{1}{2})} + {\delta_0}(b))}^2}}{nb} \frac{e^{-2nb}}{2n{b^2}} \{ 1+ O(\frac{1}{b})\}$$ i.e.
\begin{equation}\label{ub}
 {{\lVert {W_s}(nz) \rVert}_{\mathcal{P}(b, \infty)}^2} \leq 2\pi \frac{2{{(2\pi e^{(s-\frac{1}{2})} + {\delta_0}(b))}^2}}{nb} \frac{e^{-2nb}}{n{b^2}} \{ 1+ O(\frac{1}{b})\}
\end{equation}
In the last inequality, we used the following estimate from \cite[Section 3.2]{Le}: $$\displaystyle {\int_{t_1}^{t_2}} \frac{e^{-2\alpha y}}{y^2} dy = \frac{e^{-2 \alpha {t_1}}}{2 \alpha {t_1}^2} \{ 1 + O(e^{2({t_1}-{t_2})} + {t_1}^{-1}) \}$$ with an absolute constant for the $O$-term for $\alpha >1$.

Comparing \eqref{lb} and \eqref{ub} we get, for any $n \in {\mathbb Z}^*$
\begin{equation}
{{\lVert {W_s}(nz) \rVert}_{\mathcal{P}(b, \infty)}} \leq K(b) {{\lVert {W_s}(nz) \rVert}_{\mathcal{P}(2\pi, b)}}
 \end{equation}
 where $${K^2}(b)= {\frac{e^{2s-1}}{8}}{(2\pi e^{(s - \frac{1}{2})} + {\delta_0}(b))^2} e^{-|n|b}\frac{(1+ O(\frac{1}{b}))}{1+ O( e^{-|n|b} + \frac{2}{b})}.$$ From the expression it is clear that $K$ is bounded independent of $n, b$ (once $b$ is large enough) and $s \in [\frac{1}{2}, 1]$. So we obtain the claim by choosing some $b > \frac{2}{\gamma_0}$ sufficiently large (once and for all) such that the $O$-terms in the expression of $T$ are small enough. It is also clear from the expression that when $b \rightarrow \infty$, $K(b) \rightarrow 0$. This proves the Claim \ref{cusp3} and hence the first part of Lemma \ref{cusp}.

Now we prove the second part. Let $\lambda < \frac{1}{4} - \eta$ for some $\eta > 0$ and let $(\lambda, f)$ be a residual eigenpair. The Fourier expansion of $f$ inside $\mathcal{P}^1$ has the form
\begin{equation}
 f(z) = {f_0}{y^s} + \sum_{n \in {{\mathbb Z}^*}} {f_n} {W_s}(nz)= {f_0}{y^s} + g(z)
\end{equation}
where $s(1-s)=\lambda$, $s \in (0, \frac{1}{2})$ (see \cite{I}) and $g(z) = {\sum_{n \in {{\mathbb Z}^*}}} {f_n} {W_s}(nz)$. Since ${f_0}{y^s}$ and $g$ are orthogonal and since the first part can be applied to $g$, one needs only to prove the lemma for the term ${f_0}{y^s}$. So we calculate: $${\int_a^c}{y^{2s}} \frac{dy}{y^2} = \frac{1}{1-2s}\bigg(\frac{1}{{a^{1-2s}}} -\frac{1}{{c^{1-2s}}}\bigg).$$ Therefore, for $b > 2\pi$,
\begin{equation}\label{res}
{{\lVert {f_0} {y^s} \rVert}_{\mathcal{P}(b, \infty)}^2} = \frac{1}{{(\frac{b}{2\pi})^{1-2s}}-1}{{\lVert {f_0} {y^s} \rVert}_{\mathcal{P}(2\pi, b)}^2}.
\end{equation}
The lemma is satisfied by ${T_2}(b, \eta)$ such that $${T_2^2}(b, \eta) = \textrm{max}~ \bigg( {K^2}(b), \frac{1}{{(\frac{b}{2\pi})^{1-2s}} -1} \bigg).$$ From the expression it is clear that ${T_2}(b, \eta)$ depends only on two quantities: $b$ and $\frac{1}{2} - s$. Since $\frac{1}{2} - s> \sqrt{\eta} >0$, $\frac{1}{{(\frac{b}{2\pi})^{1-2s}}-1} \rightarrow 0$ when $b\rightarrow \infty$. This proves the second part.
\subsection{Mass distribution over Margulis tubes}
Now we study the distribution of the mass of a small eigenfunction over Margulis tubes. Let $\gamma$ be a simple closed geodesic of length ${l_\gamma} = 2\pi l$. Recall that ${\mathcal{C}^a}$ denotes the collar around $\gamma$ bounded by two equidistant curves of length $a$. Any $f \in {L^2}({\mathcal{C}^1})$ can be written as a Fourier series in the $\theta$-coordinate:
\begin{equation}
f(r, \theta) = {a_0}(r) + {\sum_{j=1}^{\infty}} \bigg({a_j}(r) \cos{j \theta} + {b_j}(r) \sin{j \theta}\bigg).
\end{equation}
The functions ${a_j}={a_j}(r)$ and ${b_j}={b_j}(r)$ are defined on $[-{\cosh^{-1}}(\frac{1}{l_\gamma}), {\cosh^{-1}}(\frac{1}{l_\gamma})]$ and are called the {\it $j$-th Fourier coefficients of $f$} (in $\mathcal{C}^1$). When $f$ is a $\lambda$-eigenfunction, $a_j$ and $b_j$ are solutions of the differential equation
\begin{equation}\label{uch}
\frac{{d^2}\phi}{dr^2} + {\tanh r}\frac{d\phi}{dr} + ({\lambda} - \frac{j^2}{{l^2}{{\cosh^2} r}})\phi =0.
\end{equation}
We set ${[f]_0}={a_0}(r)$ and ${[f]_1}= f - {[f]_0}$. The following lemma concerns the distribution of masses of ${[f]_0}$ and ${[f]_1}$ inside $\mathcal{C}^1$.
\begin{lem}\label{t}
For any ${l_\gamma} < \epsilon \leq {\epsilon_0}$ there exist constants ${T_1}(\epsilon), {T_2}(\epsilon) < \infty$, depending only on $\epsilon$, such that for any small eigenpair $(\lambda, f)$ of $\mathcal{C}^1$ the following inequalities hold:
\begin{equation}\label{case1}
 {{\lVert{[f]_1} \rVert}_{\mathcal{C}^\epsilon}} < {T_1}(\epsilon) {{\lVert{[f]_1}\rVert}_{{\mathcal{C}^1} \setminus {{\mathcal{C}^\epsilon}}}}
\end{equation}
and
\begin{equation}
 {{\lVert{[f]_0} \rVert}_{{\mathcal{C}^{\epsilon_0}} \setminus {\mathcal{C}^\epsilon}}} < {T_2}(\epsilon) {{\lVert{[f]_0}\rVert}_{{\mathcal{C}^1} \setminus {{\mathcal{C}^{\epsilon_0}}}}}.
\end{equation}
Therefore, for any ${l_\gamma} < \epsilon \leq {\epsilon_0}$ and any small eigenpair $(\lambda, f)$ of $\mathcal{C}^1$ one has
\begin{equation}\label{case2}
{{\lVert{f} \rVert}_{{\mathcal{C}^{\epsilon_0}} \setminus {\mathcal{C}^\epsilon}}} < \textrm{max} ~ \{{T_1}({\epsilon_0}), {T_2}(\epsilon) \} {{\lVert{f}\rVert}_{{\mathcal{C}^1} \setminus {{\mathcal{C}^{\epsilon_0}}}}}.
\end{equation}
If $\lambda < \frac{1}{4} - \eta$ for some $\eta >0$ then there exists a constant ${T_0}(\epsilon, \eta)< \infty$, depending only on $\eta$ and $\epsilon$, such that
\begin{equation}\label{case3}
 {{\lVert{[f]_0} \rVert}_{\mathcal{C}^\epsilon}} < {T_0}(\epsilon, \eta) {{\lVert{[f]_0}\rVert}_{{\mathcal{C}^1} \setminus {{\mathcal{C}^\epsilon}}}}.
\end{equation}
Furthermore, ${T_1}(\epsilon), {T_0}(\epsilon, \eta) \rightarrow 0$ as $\epsilon \rightarrow 0$.
\end{lem}
Before starting the proof of the above lemma we make a few observations about the solutions of \eqref{uch}. The change of variable $u(r) = {\cosh^{\frac{1}{2}}}(r) \phi (r)$ transforms \eqref{uch} into
\begin{equation}\label{ch}
\frac{{d^2}u}{dr^2} = \bigg((\frac{1}{4} - \lambda) + \frac{1}{4{\cosh^2}r} + \frac{j^2}{{l^2}{\cosh^2}r} \bigg) u.
\end{equation}
Let $s_j$ (resp. $c_j$) be the solution of \eqref{ch} satisfying the conditions: ${s_j}(0) = 0$ and ${s_j^{'}}(0) =1$ (resp. ${c_j}(0) = 1$ and ${c_j^{'}}(0) =0$). Since \eqref{ch} is invariant under $r \rightarrow -r$ one has: ${s_j}(-r) = -{s_j}(r)$ and ${c_j}(-r) = {c_j}(r)$ for all $j\geq 0$. Therefore there exists $t>0$ such that ${s_j} >0$ and ${c_j^{'}}>0$ on $(0, t]$. Now we prove the following claim.
\begin{claim}\label{g}
 Let $L>0$. Let $g:[0, L] \rightarrow \mathbb{R}$ be a $C^2$-function which satisfies the inequality: $$\frac{{d^2}g}{d{r^2}} > {\delta^2} g$$ for some $\delta>0$. If ${g^{'}}(0) \geq 0$ then $\displaystyle \frac{g(r)}{\cosh{\delta r}}$ is a monotone increasing function of $r$ in $(0, L]$.
\end{claim}
\textbf{Proof.} Observe that $${\bigg(\displaystyle \frac{g(r)}{\cosh{\delta r}}\bigg)^{'}} = \frac{{g^{'}}(r)\cosh{\delta r}- \delta g(r)\sinh{\delta r}}{{\cosh^2}(\delta r)}.$$ Consider the function $H$ defined on $[0, L]$ by $$H(r) = {g^{'}}(r)\cosh{\delta r}- \delta g(r)\sinh{\delta r}.$$ Since $g$ is a $C^2$ function $H$ is continuous on $[0, L]$. Observe that the claim follows if $H(r) > 0$ in $(0, L]$. Now for any $r \in (0, L]$ $${H^{'}}(r) = {g^{''}}(r)\cosh{\delta r}- {\delta^2} g(r)\cosh{\delta r}= ({g^{''}}(r) - {\delta^2} g(r)) \cosh{\delta r} > 0.$$ Therefore for $r > 0$, $H(r) > H(0) = {g^{'}}(0) \geq 0$. Hence the claim.$\Box$

\textbf{Proof of Lemma \ref{t}}. We need to estimate, for ${l_\gamma} \leq t < w \leq 1$, the quantities: $${{\lVert {{[f]_1}} \rVert}_{{\mathcal{C}^w} \setminus {\mathcal{C}^t}}^2} = {l_\gamma} {\int_{-{L_w}}^{-{L_t}}} \bigg( {\sum_{j= 1}^\infty}{\alpha_j^2} + {\beta_j^2} \bigg) dr + {l_\gamma} {\int_{{L_t}}^{{L_w}}} \bigg( {\sum_{j= 1}^\infty}{\alpha_j^2} + {\beta_j^2} \bigg) dr$$ and $${{\lVert {{[f]_0}} \rVert}_{{\mathcal{C}^w} \setminus {\mathcal{C}^t}}^2} = {l_\gamma} {\int_{-{L_w}}^{-{L_t}}} {\alpha_0^2} dr + {l_\gamma} {\int_{{L_t}}^{{L_w}}} {\alpha_0^2} dr$$ where ${\alpha_0}(r) = {\cosh^{\frac{1}{2}}}(r) {a_0}(r)$, ${\alpha_j}(r)= {a_j}(r) {\cosh^{\frac{1}{2}}}(r)$, ${\beta_j}(r)= {b_j}(r) {\cosh^{\frac{1}{2}}}(r)$ and ${L_u} = {\cosh^{-1}}(\frac{u}{l_\gamma})$. Since $s_j$ is odd and $c_j$ is even, for any symmetric subset $U \subset [-{L_1}, {L_1}]$, ${s_j}$ and ${c_j}$ are orthogonal in ${L^2}(U)$. Now ${\alpha_j}$ and ${\beta_j}$ are linear combinations of ${s_j}$ and ${c_j}$ for $j \geq 1$ and ${\alpha_0}$ is a linear combination of ${s_0}$ and ${c_0}$. Therefore, since $s_j$ and $c_j$ are orthogonal, it is enough to prove the lemma with ${s_j}$ and ${c_j}$ instead of ${[f]_1}$ and with ${s_0}$ and ${c_0}$ instead of ${[f]_0}$. We detail the computations for ${s_j}$. The computations for ${c_j}$ are similar. Let us choose $\epsilon$ such that ${l_\gamma} < \epsilon < {\epsilon_0}$. The lemma reduces to find ${K_1}(\epsilon), {K_2}(\epsilon) < \infty$, depending on $\epsilon$, and ${K_0}(\epsilon, \eta) < \infty$, depending on $\epsilon, \eta$ ($>0$), such that $${{\lVert {s_j} \rVert}_{\mathcal{C}^\epsilon}} < {K_1}(\epsilon) {{\lVert{s_j}\rVert}_{{\mathcal{C}^1} \setminus {{\mathcal{C}^\epsilon}}}}, ~ ~ ~ ~ {{\lVert{s_0} \rVert}_{{\mathcal{C}^{\epsilon_0}} \setminus {\mathcal{C}^\epsilon}}} < {K_2}(\epsilon) {{\lVert{s_0}\rVert}_{{\mathcal{C}^1} \setminus {{\mathcal{C}^{\epsilon_0}}}}}$$ and $${{\lVert{s_0} \rVert}_{\mathcal{C}^\epsilon}} < {K_0}(\epsilon, \eta) {{\lVert{s_0} \rVert}_{{\mathcal{C}^1} \setminus {{\mathcal{C}^\epsilon}}}}.$$

Let $\eta < \frac{1}{4} - \lambda$ and set ${\delta_0}= \sqrt{\eta}$ and set for $j \geq 1$, ${\delta_j}=1$. Notice that $l \cosh{r} < 1$ on $[0, {L_1})$. Hence by \eqref{ch} ${s_j}: [0, {L_1}) \rightarrow \mathbb{R}$ satisfies the inequality: $$\frac{{d^2}{s_j}}{d{r^2}} > {\delta_j^2} {s_j}.$$ Hence by Claim \ref{g} ${h_j}(r)= \frac{{s_j}(r)}{\cosh r}$, for $j \geq1$, is strictly increasing on $(0, {L_1})$. The same is true for ${h_0} = \frac{{s_0}(r)}{\cosh {\delta_0}r}$ (even when ${\delta_0}=0$).

We begin with the proof of the second part of the Lemma. So we assume $\eta>0$. For $0 \leq a < b$ consider the integral: $${\int_a^b}{s^2_0}(r) dr = {\int_a^b} {h_0^2}(r) {\cosh^2}({\delta_0} r) dr.$$ Since ${h_0}$ is strictly increasing we have
\begin{equation}\label{ieq}
 {h_0^2}(a) {\int_a^b}{\cosh^2}({\delta_0} r) dr <{\int_a^b}{s^2_0}(r) dr <{h_0^2}(b) {\int_a^b}{\cosh^2}({\delta_0} r) dr.
\end{equation}
Now choosing $a = 0$ and $b = {L_\epsilon}$ the last inequality in \eqref{ieq} gives
\begin{equation}
 {{\lVert{s_0} \rVert}_{\mathcal{C}^\epsilon}^2} < 2{l_\gamma} {h_0^2}({L_\epsilon}) {\int_0^{{L_\epsilon}}}{\cosh^2}({\delta_0} r) dr.
\end{equation}
Next choosing $a= {L_\epsilon}$ and $b={L_1}$ the first inequality in \eqref{ieq} gives \begin{equation}
{{\lVert{s_0} \rVert}_{{\mathcal{C}^1} \setminus {{\mathcal{C}^\epsilon}}}^2} > 2{l_\gamma} {h_0^2}({L_\epsilon}) {\int_{{L_\epsilon}}^{{L_1}}}{\cosh^2}({\delta_0} r) dr.
\end{equation}
Therefore
\begin{equation}\label{s0}
{{\lVert{s_0} \rVert}_{\mathcal{C}^\epsilon}} < {T_0} {{\lVert{s_0} \rVert}_{{\mathcal{C}^1} \setminus {{\mathcal{C}^\epsilon}}}}
\end{equation}
where
\begin{equation}\label{k0}
{T_0^2} = \frac{\sinh{2{\delta_0} {L_\epsilon}} + 2{\delta_0} {L_\epsilon}}{\sinh{2{\delta_0} {L_1}} - \sinh{2{\delta_0} {L_\epsilon}} + 2{\delta_0} ({L_1} - {L_\epsilon})}.
\end{equation}
We see that ${T_0}$ depends only on $\epsilon, {\delta_0}$ and $l_\gamma$. Now ${L_\epsilon}= {\cosh^{-1}}(\frac{\epsilon}{l_\gamma}) = \log (\frac{\epsilon}{l_\gamma} + \sqrt{{{(\frac{\epsilon}{l_\gamma})}^2} - 1})$. Therefore, for $\epsilon$ and ${\delta_0^2} = \eta >0 $ fixed, and $l_\gamma$ small $${T_0^2} < {K_0} \displaystyle \frac{1}{{\epsilon^{-2{\delta_0}}}-1},$$ and the constant $K_0$ is independent of $l_\gamma$ as soon as $l_\gamma$ is small compared to $\epsilon$. Thus we can choose ${T_0}(\epsilon, \eta)$ independent of $l_\gamma$ satisfying \eqref{s0}. This proves (\ref{case3})

For ${s_j}$, $j \geq 1$, exactly the same computations for $s_0$ work with ${\delta_0}$ replaced by ${\delta_j} = 1$. Hence in this case our constant, $${T_1^2}(\epsilon) < {K_1} \displaystyle \frac{1}{{\epsilon^{-2}}-1},$$ depends only on $\epsilon$. This proves \eqref{case1}.

Now we prove \eqref{case2}. Since ${s_0}: [0, {L_1}] \rightarrow {\mathbb{R}^{+}}$ is strictly increasing we have: $${\int^{L_{\epsilon_0}}_{L_\epsilon}} {s_0^2}(r) dr < {s_0^2}({L_{\epsilon_0}})({L_{\epsilon_0}} - {L_\epsilon}) ~ \textrm{and} ~ {\int_{L_{\epsilon_0}}^{L_1}} {s_0^2}(r) dr > {s_0^2}({L_{\epsilon_0}})({L_1} - {L_{\epsilon_0}}).$$ Combining the two inequalities we obtain
\begin{equation}\label{s2}
 {{\lVert{s_0} \rVert}_{{\mathcal{C}^{\epsilon_0}} \setminus {\mathcal{C}^\epsilon}}} < {T_2}(\epsilon){{\lVert{s_0} \rVert}_{{\mathcal{C}^1} \setminus {\mathcal{C}^{\epsilon_0}}}}
\end{equation}
where
\begin{equation}
 {T_2^2}(\epsilon) = \frac{{L_{\epsilon_0}} - {L_\epsilon}}{{L_1} - {L_{\epsilon_0}}}< {K_2} \bigg( \frac{\log{\frac{1}{\epsilon}}}{\log{\frac{1}{\epsilon_0}}} -1 \bigg).
\end{equation}
The constant $K_2$ is independent of $l_\gamma$ as soon as $l_\gamma$ is small compared to $\epsilon$. Thus we can choose ${T_2}(\epsilon)$ independent of $l_\gamma$ satisfying \eqref{s2}. This proves \eqref{case2}.
\subsection{Applications}
Let $S$ be a finite area hyperbolic surface with $n$ punctures. Denote by $\mathcal{P}_i$ the standard cusp around the $i$-th puncture. Recall that $\mathcal{P}_i$'s have disjoint interiors and that each of them is isometric to the half-infinite annulus $\mathcal{P}^1$ (see \ref{cuspd}). Applying Lemma \ref{cusp} in each $\mathcal{P}_i$ separately we obtain the following corollary which will be useful in our analysis.
\begin{cor}\label{crl}
For any $ 0 < \epsilon < {\epsilon_0}$ there exists $T(\epsilon) < \infty$, depending only on $\epsilon$, such that for any small cuspidal eigenpair $(\lambda, f)$ of $S$ one has
\begin{equation}\label{c1}
{{\lVert f \rVert}_{S_c^{(0, {\epsilon})}}} < T(\epsilon) {{\lVert f \rVert}_{{S_c^{(0, 1]}} \setminus {S_c^{(0, {\epsilon}]}}}}.
\end{equation}
If $\lambda < \frac{1}{4}- \eta$ for some $\eta>0$ then for any $ 0 < \epsilon < {\epsilon_0}$ there exists ${T_1}(\epsilon, \eta) < \infty$, depending only on $\epsilon$ and $\eta$, such that for any $\lambda$-eigenfunction $f$ of $S$ one has
\begin{equation}\label{c2}
{{\lVert f \rVert}_{S_c^{(0, {\epsilon})}}} < {T_1}(\epsilon, \eta) {{\lVert f \rVert}_{{S_c^{(0, 1]}} \setminus {S_c^{(0, {\epsilon}]}}}}.
\end{equation}
Furthermore, $T(\epsilon)$ and ${T_1}(\epsilon, \eta)$ tends to zero as $\epsilon \rightarrow 0$.
\end{cor}
Using this corollary and \eqref{case2} we deduce the following
\begin{cor}\label{crl2}
 For any $0 < \epsilon < {\epsilon_0}$ there exists a constant $L(\epsilon) < \infty$, depending only on $\epsilon$, such that for any small cuspidal eigenfuction $f$ of $S$ one has
 \begin{equation}
  {{\lVert f \rVert}_{S^{[\epsilon, \infty)}}} < L(\epsilon){{\lVert f \rVert}_{S^{[{\epsilon_0}, \infty)}}}.
 \end{equation}
\end{cor}
Now we give a new proof of the following theorem of D. Hejhal \cite{H}.
\begin{thm}\label{H}
Consider a sequence $({S_m}) \in {\mathcal{M}_{g, n}}$ converging to ${S_\infty} \in \overline{\mathcal{M}_{g, n}}$. Let $({\lambda_m}, {\phi_m})$ be a normalized small eigenpair of $S_m$ such that ${\lambda_m} \rightarrow {\lambda_\infty}.$ If ${\lambda_\infty} < \frac{1}{4}$ then, up to extracting a subsequence, ${\phi_m}$ converges to a normalized ${\lambda_\infty}$-eigenfunction $\phi_\infty$ of $S_\infty$.
\end{thm}
D. Hejhal's proof uses {\it convergence of Green's functions} of $S_m$ to that of $S_\infty$. Our approach is more elementary and uses the above estimates on the mass distribution of eigenfunctions over thin part of surfaces.

\textbf{Proof of Theorem \ref{H}.} First we prove that, up to extracting a subsequence, ${\phi_m}$ converges to a ${\lambda_\infty}$-eigenfunction $\phi_\infty$ of $S_\infty$. By Theorem 2 (which will be proven in \S3) it is enough to prove that there exist $\epsilon, \delta > 0$ such that ${{\lVert {\phi_m} \rVert}_{S_m^{[\epsilon, \infty)}}}\geq \delta$ up to extracting a subsequence. We argue by contradiction. Suppose that for any $\epsilon>0$ the sequence ${{\lVert {\phi_m} \rVert}_{S_m^{[\epsilon, \infty)}}}\rightarrow 0$ as $m \rightarrow \infty$. Let $\eta >0$, such that ${\lambda_m} < \frac{1}{4}-\eta$ for all $ m \geq 1$. By Lemma \ref{t} we have
\begin{equation}\label{h}
 {{\lVert {{\phi_m}} \rVert}_{\mathcal{C}^\epsilon}} < \textrm{max} \{{T_0}(\epsilon, \eta), {T_1}(\epsilon)\} {{\lVert {{\phi_m}} \rVert}_{{\mathcal{C}^1} \setminus {\mathcal{C}^\epsilon}}}.
\end{equation}
Therefore from \eqref{c2} and \eqref{h} we have
\begin{equation}
 {{\lVert {\phi_m} \rVert}_{S_m^{(0, \epsilon)}}} < \textrm{max}\{{T_0}(\epsilon, \eta), {T_1}(\epsilon), {T_1}(\epsilon, \eta)\} {{\lVert {\phi_m} \rVert}_{S_m^{[\epsilon, \infty)}}}.
\end{equation}
Hence if ${{\lVert {\phi_m} \rVert}_{S_m^{[\epsilon, \infty)}}} \rightarrow 0$ as $m \rightarrow \infty$ then ${\lVert {\phi_m} \rVert} \rightarrow 0$ as $m \rightarrow \infty$. This is a contradiction to the fact that each $\phi_m$ is normalized i.e. ${\lVert {\phi_m} \rVert} = 1$.

Next we prove that ${\lVert {\phi_\infty} \rVert}= 1.$ By uniform convergence over compacta, in each cusp and in each pinching collar, the Fourier coefficients of $\phi_m$ will converge to the corresponding Fourier coefficients of $\phi_\infty$. Therefore, by \eqref{case1}, \eqref{case3} and \eqref{c2}, $\phi_m$'s are uniformly integrable: for any $\delta > 0$ there exist $\epsilon > 0$ such that for all large values of $m$
\begin{equation}
{{\lVert {\phi_m} \rVert}_{S_m^{[\epsilon, \infty)}}} > 1 - \delta.
\end{equation}
Hence ${\lVert {\phi_\infty} \rVert}= 1.$ This finishes the proof. $\Box$
\section{Proof of Theorem 2}
Let $({S_m})$ be a sequence in ${\mathcal{M}_{g, n}}$ which converges in $\overline{{\mathcal{M}_{g, n}}}$ to ${S_\infty}$. Let ${\Gamma_m}, {\Gamma_\infty}$ be such that ${S_m}= {\mathbb{H} / {\Gamma_m}}$ and ${S_\infty}= {\mathbb{H}/ {\Gamma_\infty}}$. Recall that the convergence ${S_m} \rightarrow {S_\infty}$ means that for any fixed positive constant ${\epsilon_1} \leq {\epsilon_0}$ ($\epsilon_0$ is the Margulis constant) and a choice of base point ${p_m} \in {S_m^{[{\epsilon_1}, \infty)}}$, after conjugating $\Gamma_m$ so that the projection $\mathbb{H} \rightarrow \mathbb{H}/{\Gamma_m}$ maps $i$ to $p_m$, $(\mathbb{H}/{\Gamma_m}, {p_m})$ converges to a component $(\mathbb{H}/{\Gamma_\infty}, {p_\infty})$ of ${S_\infty}$. We begin by fixing some $\epsilon < {\epsilon_0}$ and ${p_m} \in {S_m^{[\epsilon, \infty)}}$. In the following we assume that $\epsilon_1$, $p_m$, $\Gamma_m$, $p_\infty$ and $\Gamma_\infty$
satisfy the previous statement.

To simplify notations we shall assume that only one closed geodesic $\gamma_m$ gets pinched as ${S_m} \rightarrow {S_\infty} \in \partial{\mathcal{M}_{g, n}}$. In particular the limit surface $S_\infty$ (which may be disconnected) has two new cusps. Denote the standard cusps of $S_m$ by ${\mathcal{P}_1}(m),$ ${\mathcal{P}_2}(m),..., {\mathcal{P}_n}(m)$ and the limits of these in ${S_\infty} \in {\partial {\mathcal{M}_{g, n}}}$ by ${\mathcal{P}_1}(\infty)$,...,$ {\mathcal{P}_n}(\infty)$ and denote by ${\mathcal{P}_{n+ 1}}(\infty), {\mathcal{P}_{n+2}}(\infty)$ the {\it new cusps} which arise due to the pinching of $\gamma$. The cusps ${\mathcal{P}_i}(\infty)$ for $1 \leq i \leq n$ will be called {\it old cusps}.

Recall that we have a sequence of small cuspidal eigenpairs $({\lambda_m}, {\phi_m})$ of ${S_m}= {\mathbb{H} / {\Gamma_m}}$ such that the $L^2$-norm of $\phi_m$ is 1 and ${\lambda_m} \rightarrow {\lambda_\infty} \leq \frac{1}{4}$.
\begin{notn}
In what follows $d{\mu_m}$ will denote the area measure on $S_m$ for $m \in \mathbb{N}\cup\{\infty\}$ and $d\mu_{\mathbb{H}}$ will denote the area measure on $\mathbb{H}$. The lift of $f \in {L^2}({S_m})$ to $\mathbb{H}$ under the projection $\mathbb{H} \rightarrow \mathbb{H}/ {\Gamma_m}$, defined as above, will be denoted by $\widetilde{f}$.
\end{notn}
By Green's formula one has: $$\int_{{S_m}} |\nabla{\phi_m}|^2 d{\mu_m} = {\lambda_m} \int_{{S_m}} |{\phi_m}|^2 d{\mu_m} = {\lambda_m}.$$ Let $K \subset \mathbb{H}$ be compact. One can cover $K$ by finitely many geodesic balls of radius $\rho$. If $\rho$ is sufficiently small then each of these balls maps injectively to $S_m$ since ${\Gamma_m}\rightarrow {\Gamma_\infty}$. Therefore, since $\lVert {\phi_m} \rVert = 1$ $\lVert {\widetilde{\phi_m}|_K} \rVert$ is bounded depending only on $K$. From the mean value formula \cite[Corollary 1.3]{F} there exists a constant $\Lambda({\lambda_\infty}, \rho)$ such that for $\lambda_m$ close to $\lambda_\infty$, $$|\widetilde{\phi_m}(q)| \leq \Lambda({\lambda_\infty}, \rho) {\int_{N(K, \frac{\rho}{2})}} |\widetilde{\phi_m}| d{\mu_{\mathbb{H}}}$$ for each $q \in K$ where $N(K, r)$ denotes the closed neighborhood of radius $r$ of $K$ in $\mathbb{H}$. Next we use the $L^p$-Schauder estimates \cite[Theorem 4, Sect. II.5.5]{B-J-S} to obtain a uniform bound for $\nabla\widetilde{\phi_m}$ on $N(K, \frac{\rho}{2})$. This makes $({\widetilde{\phi_m}|_K})$ an equicontinuous family. So, by Arzela-Ascoli theorem, up to extracting a subsequence, $(\widetilde{\phi_m})$ converges to a continuous function $\widetilde{\phi_\infty}$ on $K$. By a diagonalization argument one may suppose that the sequence works for all compact subsets of $\mathbb{H}$. Therefore, up to extracting a subsequence, $\widetilde{\phi_m} \rightarrow \widetilde{\phi_\infty}$ uniformly over compacta. By this uniform convergence it is clear that $\widetilde{\phi_\infty}$ is a {\it weak solution} of the Laplace equation: $\Delta u + {\lambda_\infty}u=0$. Therefore, by elliptic regularity, $\widetilde{\phi_\infty}$ indeed a smooth and satisfies $${\Delta} \widetilde{\phi_\infty} + {\lambda_\infty} \widetilde{\phi_\infty} =0.$$ Also by the convergence $\widetilde{\phi_\infty}$ induces a function $\phi_\infty$ on $S_\infty$ that satisfies $$\Delta{\phi_\infty} + {\lambda_\infty}{\phi_\infty} =0.$$  However, $\phi_\infty$ may not be an eigenfunction since it could be the zero function. In order to discuss this point, we shall consider two cases according to whether the ${L^2}$-norm ${{\lVert {\phi_m} \rVert}_{S_m^{[\epsilon, \infty)}}}$ of the restriction of ${\phi_m}$ to ${S_m^{[\epsilon, \infty)}}$ is bounded below by a positive constant or not.
\subsubsection*{Case 1: $\exists$ $\epsilon, \delta > 0$ such that $\limsup{{\lVert {\phi_m} \rVert}_{S_m^{[\epsilon, \infty)}}} \geq \delta$.}
We may assume that $\lim{{\lVert {\phi_m} \rVert}_{S_m^{[\epsilon, \infty)}}} \geq \delta$. Then by the uniform convergence of $\widetilde{\phi_m} \rightarrow \widetilde{\phi_\infty}$ over compacta, $${\int_{S_\infty^{[\epsilon, \infty)}}} {\phi_\infty^2} d{\mu_\infty} = {\lim_{{m_j} \rightarrow \infty}} {\int_{S_{m_j}^{[\epsilon, \infty)}}}{\phi_{m_j}^2} d{\mu_{m_j}} \geq \delta>0.$$ Therefore $\phi_\infty$ is not the zero function and its $L^2$ norm is less than $1$. Therefore it is a ${\lambda_\infty}$-eigenfunction.
\subsubsection*{Case 2: For any $\epsilon > 0$ the sequence ${{\lVert {\phi_m}\rVert}_{S_m^{[\epsilon, \infty)}}} \rightarrow 0$.}
Then we will prove the following statements: $\\ (i)$ ${S_\infty} \in \partial {\mathcal{M}_{g, n}}$, $\\(ii)$ ${\lambda_\infty} = \frac{1}{4}$ and $\\(iii)$ $\exists$ constants $K_m$ such that, up to extracting a subsequence, $({K_m}\widetilde{\phi_m})$ converges uniformly to a function which is a linear combination of Eisenstein series and (possibly) a $\frac{1}{4}$-cuspidal eigenfunction.

\underline{(i)}~ Suppose by contradiction that ${S_\infty} \in {\mathcal{M}_{g, n}}$. Then all the cusps of $S_\infty$ are old cusps. Let $s({S_\infty})$ denote the {\it systole} of $S_\infty$. Then, for $0< \epsilon < \frac{s({S_\infty})}{2}$ and for $m$ large enough, we have ${S_m^{(0, \epsilon)}} \subset {\cup_{i=1}^n}{\mathcal{P}_i}(m)$. Therefore, applying Corollary \ref{crl}, the assumption ${{\lVert {\phi_m}\rVert}_{S_m^{[\epsilon, \infty)}}} \rightarrow 0$ implies that ${\lVert {\phi_m} \rVert} \rightarrow 0$. This is a contradiction since each $\phi_m$ is normalized. Thus ${S_\infty} \in  \partial{\mathcal{M}_{g, n}}$.

\underline{(ii)} follows from Theorem \ref{H}.

\underline{(iii)}~ Fix some $\epsilon$, $0 < \epsilon < {\epsilon_0}$. Choose constants $K_m \geq 1$ such that $${\int_{S_m^{[\epsilon, \infty)}}} |{K_m}{\phi_m}|^2 d{\mu_m} = 1.$$ Therefore the sequence $({K_m})$ must diverge to $\infty$. Using mean value formula \cite[Corollary 1.3]{F}, ${L^p}$-Schauder estimates \cite{B-J-S} and elliptic regularity, as earlier, and Corollary \ref{crl2} we obtain that, up to extracting a subsequence, $(\widetilde{{K_m}{\phi_m}})$ converges, uniformly over compacta, to a $C^\infty$ function $\widetilde{\phi_\infty}$ that satisfies $$\Delta {\widetilde{\phi_\infty}} + \frac{1}{4} {\widetilde{\phi_\infty}} =0.$$ Moreover, $\widetilde{\phi_\infty}$ induces a function $\phi_\infty$ on $S_\infty$ that satisfies
\begin{equation}
{\Delta} {\phi_\infty} + \frac{1}{4} {\phi_\infty} = 0.
\end{equation}
Using the uniform convergence over compacta we have $${\int_{S_\infty^{[\epsilon, \infty)}}} {\phi_\infty^2} d{\mu_\infty} = {\lim_{m \rightarrow \infty}} {\int_{S_m^{[\epsilon, \infty)}}}{K_m}{\phi_m^2} d{\mu_m} = 1.$$ Therefore $\phi_\infty$ is not the zero function. From Lemma \ref{cusp} and Lemma \ref{t} \eqref{case1} we deduce that ${\phi_\infty}$ satisfies {\it moderate growth} condition \cite[p. 80]{Wo} in each cusp. It is known that for any $\lambda \geq \frac{1}{4}$ the space of {\it moderate growth $\lambda$-eigenfunctions} of $S_\infty$ is spanned by Eisenstein series and (possibly) $\lambda$-cuspidal eigenfunctions (see \S3 in \cite{Wo}). In particular, $\phi_\infty$ is a linear combination of Eisenstein series and (possibly) a cuspidal eigenfunction. This finishes the proof of $(iii). \Box$
\section{Proof of Theorem 1}
We begin by proving \underline{Lemma 1} which says that ${\mathcal{C}_{g, n}^{\frac{1}{4}}}(k)$ is open in $\mathcal{M}_{g, n}$.
\subsection{Proof of \underline{Lemma 1}}
Empty set is open by convention. Therefore, we argue by contradiction and assume that there exists a $S \in {\mathcal{C}_{g, n}^{\frac{1}{4}}}(k)$ such that every neighborhood of $S$ contains points from ${\mathcal{M}_{g, n}} \setminus {\mathcal{C}_{g, n}^{\frac{1}{4}}}(k)$. In other words, there exists a sequence $({S_m}) \subseteq {\mathcal M}_{g, n}$ that converges to $S$ and, for all $m$, ${\lambda^c_{k}}(S_m) \leq \frac{1}{4}$. For $1 \leq i \leq k$, let us denote by $\phi^i_m$ a normalized ${\lambda^c_i}(S_m)$-cuspidal eigenfunction such that ${\{{\phi^i_m}\}_{i=1}^{k}}$ is an orthonormal family in ${L^2}(S_m)$. Since we are considering small eigenvalues, up to extracting a subsequence, the sequence $({\lambda^c_i}(S_m))$ converges. For simplicity we assume that, for $1 \leq i \leq k$, the sequence $({\lambda^c_i}(S_m))$ converges and denote by $\lambda^i_\infty$ its limit. Observe that, for $1 \leq i \leq k$, ${\lambda^i_\infty} \leq \frac{1}{4}.$ Now, since $S \in {\mathcal M}_{g, n}$ by Theorem 2, up to extracting a subsequence, $({\phi^i_m})$ converge to $\lambda^i_\infty$-eigenfunction $\phi^i_\infty$ of $S$. Moreover, by the result about uniform integrability inside cusps in Corollary \ref{crl}: $\lVert {\phi_\infty^i} \rVert = 1$. Hence ${\{{\phi^i_\infty}\}_{i=1}^{k}}$ is an orthonormal family in ${L^2}(S)$ so that the $k$-th cuspidal eigenvalue ${\lambda^c_{k}}(S)$ of $S$ is below $\frac{1}{4}$. This is a contradiction because by our assumption ${\lambda^c_{k}}(S) > \frac{1}{4}$ as $S \in {\mathcal{C}_{g, n}^{\frac{1}{4}}}(k)$.$\Box$

Now we give a proof of Proposition \ref{neighbrhd} which says

\textbf{Proposition \ref{neighbrhd}}~ {\it There exists a neighborhood $\mathcal{N}({\mathcal{M}_{g, 1}} \cup {\mathcal{M}_{0, n+1}})$ of ${\mathcal{M}_{g, 1}} \cup {\mathcal{M}_{0, n+1}}$ in $\overline{\mathcal{M}_{g, n}}$ such that for each $S \in \mathcal{N}({\mathcal{M}_{g, 1}} \cup {\mathcal{M}_{0, n+1}})$: ${\lambda^c_{2g-1}}(S) > \frac{1}{4}$ i.e. $$\mathcal{N}({\mathcal{M}_{g, 1}} \cup {\mathcal{M}_{0, n+1}}) \subset {{\mathcal{C}_{g, n}^{\frac{1}{4}}}}(2g-1).$$}
\subsection{Proof of Proposition}
We argue by contradiction and assume that there is a sequence ${S_m} \in {\mathcal{M}_{g, n}}$ converging to ${S_\infty} \in {\mathcal{M}_{g, 1}} \cup {\mathcal{M}_{0, n+1}} \subset \partial{\mathcal{M}_{g, n}}$ such that ${\lambda_{2g-1}^c}({S_m}) \leq \frac{1}{4}$. For $1 \leq i \leq 2g-1$ and for each $m$ we choose small cuspidal eigenpairs $({\lambda^i_m}, {\phi^i_m})$ of $S_m$ such that $\\(i)~ \{{\phi^i_m}\}_{i=1}^{2g-1}$ is an orthonormal family in ${L^2}({S_m}),$ $\\(ii)~ {\lambda_m^i}$ is the $i$-th eigenvalue of $S_m$.

Theorem 2 provides two possible behaviors of the sequence $({\phi_m^i})$. However in our case we have Lemma 2:

\underline{\textbf{Lemma 2}}~ {\it For each $i$, $1 \leq i \leq 2g-1$, up to extracting a subsequence, the sequence $({\phi^i_m})$ converges to a $\lambda^i_\infty$-eigenfunction $\phi^i_\infty$ of $S_\infty$. The limit functions ${\phi^i_\infty}$ and ${\phi^j_\infty}$ are orthogonal for $i \neq j$ i.e. $S_\infty$ has at least $2g-1$ small eigenvalues. Moreover none of the ${\phi^i_\infty}$ is residual.}
\subsubsection{Proof of Lemma 2}\label{pc1}
By uniform convergence of $\phi^i_m$ to $\phi^i_\infty$, we have $\lVert {\phi^i_\infty} \rVert \leq 1$. To prove the first two statements of the lemma it is enough to prove that, for $1 \leq i \leq 2g-1$, $\lVert {\phi^i_\infty} \rVert = 1$ because this will imply that $\phi_\infty^i$ is not the zero function and that $({\phi^i_m})$ is uniformly integrable over the thick parts: for any $t >0$ there exists $\epsilon$ such that for all $m$ one has, $${{\lVert {\phi_m} \rVert}_{S_m^{[\epsilon, \infty)}}} > 1 - t.$$ To prove that, for each $1 \leq i \leq 2g-1$, $\lVert {\phi^i_\infty} \rVert = 1$ we argue by contradiction and assume that for some $1 \leq i \leq 2g-1$, $\lVert {\phi^i_\infty} \rVert = 1 - \delta.$ To simplify the notation, denote the sequence $({\lambda^i_m}, {\phi^i_m})$ by $({\lambda_m}, {\phi_m})$ and the limit $({\lambda^i_\infty}, {\phi^i_\infty})$ by $({\lambda_\infty}, {\phi_\infty})$. By Corollary \ref{crl} the functions $\phi_m$ are uniformly integrable over the union of cusps of $S_m$: for any $t > 0$ there exists $\epsilon > 0$ such that for all $m$ one has:
\begin{equation}\label{uniform}
{{\lVert {\phi_m} \rVert}_{S_{m, c}^{(0, \epsilon)}}}< t.
\end{equation}
Since ${S_\infty} \in {\mathcal{M}_{g, 1}} \cup {\mathcal{M}_{0, n+1}}$ there is only one closed geodesic, ${\gamma_m} \subset {S_m}$, whose length $l_{\gamma_m}$ tends to zero. For any $l \leq 1$ and for $m$ large enough such that ${l_{\gamma_m}} < l $ denote by ${\mathcal{C}^l_ m} \subset {S_m}$ the collar around $\gamma_m$ bounded by two equidistant curves of length $l$. In view of the uniform integrability inside cusps \eqref{uniform}, there exists ${\epsilon_0}>0$ such that for any $\epsilon \leq {\epsilon_0}$ there exists $m(\epsilon)$ such that for $m \geq m(\epsilon)$ we have:
\begin{equation}\label{discoll}
{{\lVert {\phi_m} \rVert}_{\mathcal{C}^\epsilon_m}} > \frac{\delta}{2}.
\end{equation}
Now we distinguish again two cases depending on whether ${\lambda_\infty}< \frac{1}{4}$ or ${\lambda_\infty} = \frac{1}{4}$. If ${\lambda_\infty} < \frac{1}{4}$ then we have a contradiction since $\lVert {\phi_\infty} \rVert = 1$ by Theorem \ref{H}. Hence we may suppose that ${\lambda_\infty} = \frac{1}{4}$. So, by Theorem 2 either $\phi_\infty$ is the zero function or, for instance by \cite[Theorem 3.2]{I}, $\phi_\infty$ is cuspidal. Now recall that by lemma \ref{t} we have uniform integrability of $[{\phi_m}]_1$: for any $t$ there exists $\epsilon$ such that for all $m$: $${{\lVert {[{\phi_m}]_1} \rVert}_{\mathcal{C}^\epsilon_m}} < t.$$ Hence by \eqref{discoll}, there exists ${\epsilon_1}$ such that for any $\epsilon \leq {\epsilon_1}$ the exists ${m_1}(\epsilon)$ such that for $m \geq {m_1}(\epsilon)$ one has:
\begin{equation}\label{phi0}
{{\lVert {[{\phi_m}]_0} \rVert}_{\mathcal{C}^\epsilon_m}} > \frac{\delta}{4}
\end{equation}
In particular, if $c(\epsilon, m) = {\textrm{sup}_{z \in {\mathcal{C}^\epsilon_m}}} |{[{\phi_m}]_0}|$ then, since area of $\mathcal{C}^\epsilon_m$ is less than $1$, we have for any $\epsilon \leq {\epsilon_1}$ and $m \geq {m_1}(\epsilon)$:
\begin{equation}\label{max}
c(\epsilon, m) > \frac{\delta}{4}.
\end{equation}
Now we prove that $[{\phi_m}]_1$ is uniformly small inside $\mathcal{C}^\epsilon_m$. More precisely,
\begin{lem}\label{pcurve}
Let $\epsilon$ be such that $0< \epsilon <1$. There exists a constant $K < \infty$, independent of $\epsilon$, and ${m_2}(\epsilon) \in \mathbb{N}$ such that for $m \geq {m_2}(\epsilon)$ and $z \in {\mathcal{C}^\epsilon_m}$: $$|{[{\phi_m}]_1}|(z) < K\frac{\epsilon^\frac{1}{2}}{1 - \epsilon}.$$
\end{lem}
\textbf{Proof.} Consider the expansion of ${\phi_m}$ inside ${\mathcal{C}^1_m}$ with respect to the Fermi coordinates (see \ref{cylinder}):
\begin{equation}\label{expr}
{\phi_m}(r, \theta) = {a_0^m}(r) + {\sum_{j=1}^{\infty}} \bigg({a_j^m}(r) \cos{j \theta} + {b_j^m}(r) \sin{j \theta}\bigg).
\end{equation}
Here, for each $j \geq 0$, $({a_j^m}, {b_j^m})$ are the $j$-th Fourier coefficients of $\phi_m$ inside $\mathcal{C}^1_ m$ and are defined for all $|r| \leq {L_{1, m}}$. Recall that, for any $\epsilon \in [{l_{\gamma_m}}, 1]$ we denote by ${L_{\epsilon, m}}$ the number ${\cosh^{-1}}(\frac{\epsilon}{l_{\gamma_m}})$. Recall also that since ${\phi_m}$ is a ${\lambda_m}$-eigenfunction, ${a_j^m}$ and ${b_j^m}$ satisfy \eqref{uch} with $2\pi l= {l_{\gamma_m}}$ and $\lambda= {\lambda_m}$. Therefore, for $j \geq 1$, one can express:
\begin{equation}\label{expression}
(1) ~ {a_j^m}(r) = {a_{m, j}}{s_{m, j}}(r) + {b_{m, j}}{c_{m, j}}(r) $$$$
(2) ~ {b_j^m}(r) = {{a_{m, j}}^{'}}{s_{m, j}}(r) + {{b_{m, j}}^{'}}{c_{m, j}}(r)
\end{equation}
where ${s_{m, j}}(r)$ and ${c_{m, j}}(r)$ are the two linearly independent solutions of \eqref{uch} with $l = l({\gamma_m})$ and $\lambda = {\lambda_m}$.

Recall that ${s_{m, j}}(r){\cosh^{\frac{1}{2}}}(r)$ and ${c_{m,j}}(r){\cosh^{\frac{1}{2}}}(r)$ satisfy:
$$\frac{{d^2}u}{dr^2} = \bigg(\frac{1}{4{\cosh^2}r} + \frac{{j^2}}{{l^2}{\cosh^2}r} \bigg) u.$$
Since, for $r \leq {L_{\epsilon, m}}$, ${{l^2}{\cosh^2}r} \leq 1$ by Claim \ref{g}, for each $j \geq 1$, there exists strictly increasing functions ${h_{m, j}}: [0, {L_{1, m}}] \rightarrow {\mathbb{R}_{>0}}$ and ${k_{m, j}}: [0, {L_{1, m}}] \rightarrow {\mathbb{R}_{>0}}$ such that
\begin{equation}\label{sj}
(i) ~ {s_{m, j}}(r)\sqrt{\cosh(r)} = {h_{m, j}}(r) \cosh jr $$$$
(ii) ~ {c_{m, j}}(r)\sqrt{\cosh(r)}= {k_{m, j}}(r)\cosh jr.
\end{equation}
We denote by ${\mathcal{P}_{n+1}}(\infty)$ and ${\mathcal{P}_{n+2}}(\infty)$ the two new cusps of $S_\infty$ that appear as the limit of $\mathcal{C}^1_m$ as $m \rightarrow \infty$. Now, let us assume:$${\textrm{sup}_{z \in \partial{\mathcal{P}_{n+1}}(\infty) \cup \partial{\mathcal{P}_{n+2}}(\infty)}}|{\phi_\infty}|(z) < \frac{t}{4}.$$ Then, by the uniform convergence of ${\phi_m}$ to ${\phi_\infty}$ over compacta, we have a $N \in \mathbb{N}$ such that for $ m \geq N$ and $z \in \partial {\mathcal{C}^1_m}$: $$|{\phi_m}|(z) < \frac{t}{4}.$$ By \eqref{expr} for any $j \geq 1$:
\begin{equation}\label{am}
|{a^m_j}|(\pm {L_{1, m}}) = \frac{1}{\pi}|{\int_0^{2\pi}}{\phi_m}(\pm {L_{1, m}}, \theta) \cos{j\theta} d\theta| \leq \frac{t}{2}.
\end{equation}
Similar calculations for $b^m_j$ provide: $|{b^m_j}|(\pm {L_{1, m}}) \leq \frac{t}{2}$. Recall that ${s_{m, j}}$ is odd and ${c_{m, j}}$ is even. So by \eqref{expression} and \eqref{sj}:
\begin{equation}\label{sum}
(i) ~ {a^m_j}({L_{1, m}}) + {a^m_j}(-{L_{1, m}}) = 2 {b_{m, j}} {k_j}({L_{1, m}}) \frac{\cosh {j{L_{1, m}}}}{\sqrt{\cosh {L_{1, m}}}} $$$$
(ii) ~ {a^m_j}({L_{1, m}}) - {a^m_j}(-{L_{1, m}}) = 2 {a_{m, j}} {h_j}({L_{1, m}}) \frac{\cosh {j{L_{1, m}}}}{\sqrt{\cosh {L_{1, m}}}}.
\end{equation}
Therefore, by \eqref{am} and \eqref{sum}:
\begin{equation}\label{ineq}
(i) ~ |{b_{m, j}}|{k_j}({L_{1, m}})\frac{\cosh {j{L_{1, m}}}}{\sqrt{\cosh {L_{1, m}}}} < \frac{t}{2} $$$$
(ii) ~ |{a_{m, j}}|{h_j}({L_{1, m}})\frac{\cosh {j{L_{1, m}}}}{\sqrt{\cosh {L_{1, m}}}} < \frac{t}{2} .
\end{equation}
Therefore, for any $r \leq {L_{1, m}}$: $$|{a^m_j}|(r) =|{a_{m, j}}{s_{m, j}}(r) + {b_{m, j}}{c_{m, j}}(r)| < |{a_{m, j}}|{s_{m, j}}(r) + |{b_{m, j}}|{c_{m, j}}(r).$$ The last term of the inequality is $$|{a_{m, j}}| {h_{m, j}}(r) \frac{\cosh jr}{\sqrt{\cosh r}} + |{b_{m, j}}|{k_{m, j}}(r) \frac{\cosh jr}{\sqrt{\cosh r}} < t \frac{\cosh jr}{\sqrt{\cosh r}}\frac{\sqrt{\cosh {L_{1, m}}}}{\cosh j {L_{1, m}}}$$ since ${h_{m, j}}$ and ${k_{m, j}}$ are strictly increasing functions (by \eqref{ineq}). Similarly, $$|{b^m_j}|(r) < t \frac{\cosh jr}{\sqrt{\cosh r}}\frac{\sqrt{\cosh {L_{1, m}}}}{\cosh j {L_{1, m}}}.$$ Hence
\begin{equation}\label{infsum}
|{[{\phi_m}]_1}|(r, \theta) < 2t {\sum_{j=1}^\infty}\frac{\cosh jr}{\sqrt{\cosh r}}\frac{\sqrt{\cosh {L_{1, m}}}}{\cosh j {L_{1, m}}}.
\end{equation}
Since, for $j \geq 1$, the function $\frac{\cosh jr}{\sqrt{\cosh r}}$ is strictly increasing, for any $r \leq {L_{\epsilon, m}}$ :
\begin{equation}\label{monot}
{\sum_{j=1}^\infty}\frac{\cosh jr}{\sqrt{\cosh r}}\frac{\sqrt{\cosh {L_{1, m}}}}{\cosh j {L_{1, m}}} < {\sum_{j=1}^\infty}\frac{\cosh j{L_{\epsilon, m}}}{\sqrt{\cosh {L_{\epsilon, m}}}}\frac{\sqrt{\cosh {L_{1, m}}}}{\cosh j {L_{1, m}}}
\end{equation}
Now fix an $\epsilon$ such that $0 < \epsilon < 1$. Observe that ${L_{\epsilon, m}}= \log(\frac{\epsilon}{l_{\gamma_m}} + \sqrt{{(\frac{\epsilon}{l_{\gamma_m}})^2}-1})$. So, for $m$ large such that $l_{\gamma_m}$ is small compared to $\epsilon$:
\begin{equation}\label{cosh}
{\sum_{j=1}^\infty}\frac{\cosh j{L_{\epsilon, m}}}{\sqrt{\cosh {L_{\epsilon, m}}}}\frac{\sqrt{\cosh {L_{1, m}}}}{\cosh j {L_{1, m}}} < {K^{'}} {\sum_{j=1}^\infty} {\epsilon^j}{\epsilon^{-\frac{1}{2}}} = {K^{'}} \frac{\epsilon^{\frac{1}{2}}}{1- \epsilon}
\end{equation}
where the constant ${K^{'}}$ can be chosen independently of $\epsilon$ as soon as $m$ is larger than some number ${m_2}(\epsilon) \in \mathbb{N}$. Therefore, by \eqref{infsum} and \eqref{cosh}, for $m \geq {m_2}(\epsilon)$ and $(r, \theta) \in {\mathcal{C}^\epsilon_m}$
\begin{equation}\label{uniformb}
|{[{\phi_m}]_1}|(r, \theta) < 2t{K^{'}} \frac{\epsilon^{\frac{1}{2}}}{1- \epsilon}.
\end{equation}
This proves the lemma.$\Box$

Now fix $\epsilon < {\epsilon_1}$ (see \eqref{phi0}) such that $K \frac{\epsilon^\frac{1}{2}}{1- \epsilon} < \frac{\delta}{4}$ and choose $m \geq \textrm{max} \{ {m_1}(\epsilon), {m_2}(\epsilon) \}$. Then by Lemma \ref{pcurve} and \eqref{max}: for each $z \in {\mathcal{C}^\epsilon_m}$
\begin{equation}
c(\epsilon, m) > |{[{\phi_m}]_1}|(z).
\end{equation}
So the parallel curve $\alpha_m$ with distance ${r_0}$ ($\leq {L_{\epsilon, m}}$) from $\gamma_m$ such that $c =|{[{\phi_m}]_0}|({r_0})$ has the property that $\phi_m$ has constant sign on it. In other words, the nodal set $\mathcal{Z}({\phi_m})$ does not intersect this curve. This is a contradiction to the next lemma.
\begin{lem}\label{curve}
Let $S$ be a noncompact, finite area hyperbolic surface of type $(g, n)$. Let $\gamma$ be a simple closed geodesic that separates $S$ into two connected components $\mathcal{T}_1$ and $\mathcal{T}_2$ such that $\mathcal{T}_1$ is topologically a sphere with $n+1$ punctures and $\mathcal{T}_2$ is topologically a genus $g$ surface with one puncture. Let $f$ be a small cuspidal eigenfunction of $S$. Then the zero set ${\mathcal Z}(f)$ of $f$ intersects every curve homotopic to $\gamma$.
\end{lem}

\centerline{\includegraphics[height=3in]{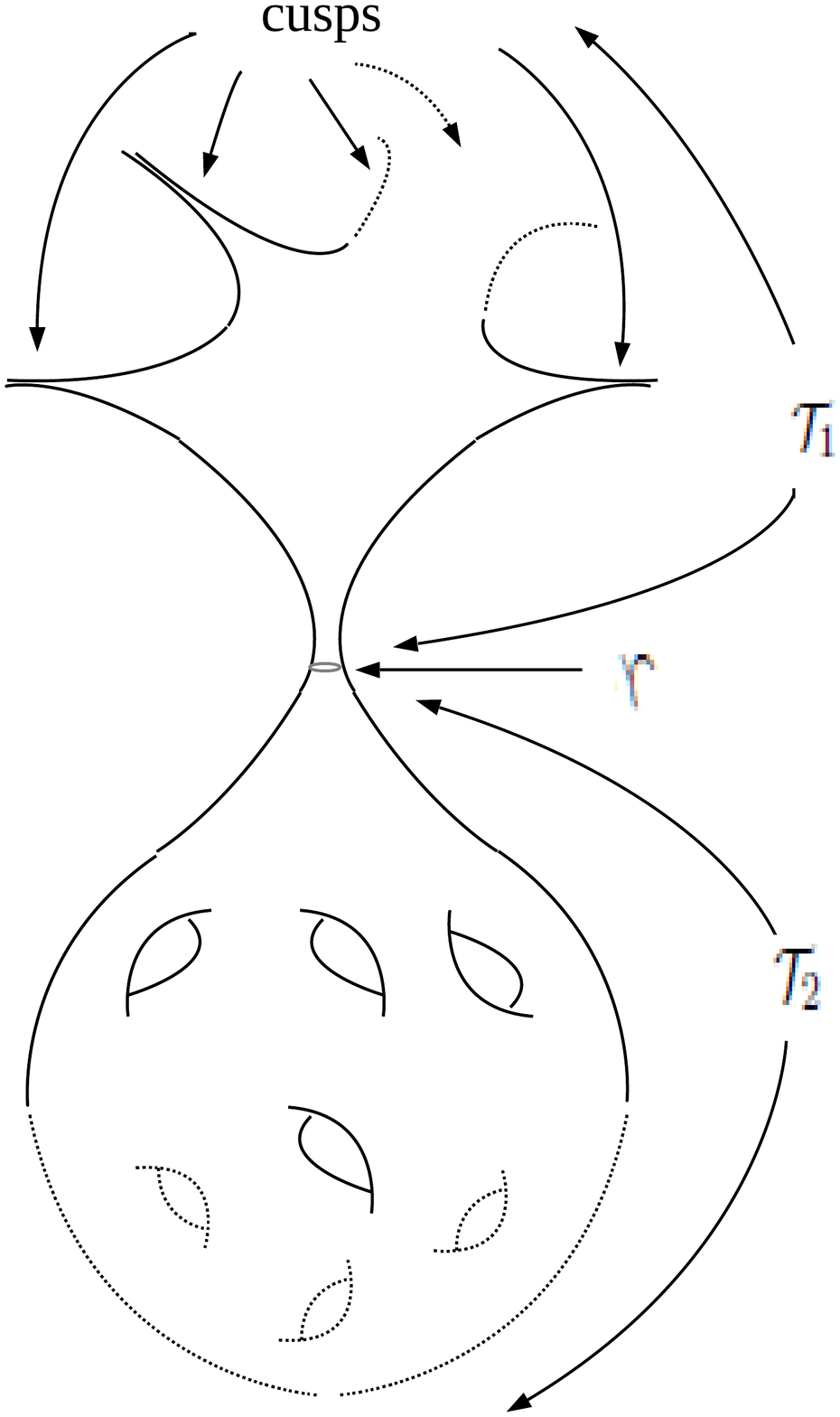}}

\textbf{Proof.} Recall that ${\mathcal Z}(f)$ is a locally finite graph \cite{Ch}. Let us assume that ${\mathcal Z}(f)$ does not intersect some curve $\tau$ homotopic to $\gamma$. We have ${S \setminus \tau} = {\mathcal{T}_1} \cup {\mathcal{T}_2}$ and all the punctures of $S$ are contained in ${\mathcal T}_1$. Consider the components of ${{\mathcal T}_1} \setminus {\mathcal Z}(f)$. Recall that since $f$ is cuspidal $\mathcal{Z}(f)$ contains all the punctures of $S$ and therefore these components give rise to a cell decomposition of a once punctured sphere. The Euler characteristic of the component $\mathcal{F}$ containing $\tau$ as a puncture is either negative or zero (since $\gamma$ and each component of $\mathcal{Z}(f)$ are essential; see \cite{O}). Each component of ${{\mathcal T}_1} \setminus {\mathcal Z}(f)$ other than $\mathcal{F}$ (at least one such exists since $g$ changes sign in $\mathcal{T}_1$) is a {\it nodal domain} of $f$ and hence has negative Euler characteristic \cite{O}. Also ${\mathcal Z}(f)$ being a graph has non-positive Euler characteristic. Let $C^{+}$ (resp. $C^{-}$) be the union of the nodal domains contained in ${\mathcal T}_1$ which are different from $\mathcal{F}$ and where $f$ is positive (resp. negative). Denote by $\chi(X)$ the Euler characteristic of the topological space $X$. Since the Euler characteristic of a once punctured sphere is $1$, by the Euler-Poincar\'{e} formula one has: $$1 = \chi(\mathcal{F})+ \chi( C^{+}) + \chi(C^{-}) + \chi({\mathcal Z}(f)).$$ This is a contradiction because the right hand side of the equality is strictly negative.$\Box$

Now we prove that $\phi_\infty$ is not a residual eigenfunction. It is clear from the uniform convergence that ${\phi_\infty}$ is cuspidal at the old cusps. If ${\phi_\infty}$ is a residual eigenfunction then the only possibility is that $\phi_\infty$ is not cupsidal at one of the two new cusps. Let us assume that ${\phi_\infty}$ is residual in $\mathcal{P}_{n+1}$. Then, for sufficiently large $t$, ${\phi_\infty}$ has constant sign in $\mathcal{P}_{n+1}^t$. Therefore, by the uniform convergence ${{\phi_m}}|_{S_m^{[\epsilon, \infty)}} \rightarrow {{\phi_\infty}|_{S_\infty^{[\epsilon, \infty)}}}$ it follows that, for all $m$ large, ${\phi_m}$ has constant sign on a component of $\partial{\mathcal{C}_m^\frac{1}{t}}$. Since this component is homotopic to $\gamma_m$ this leads to a contradiction to Lemma \ref{curve} as well. This finishes the proof of Lemma 2.$\Box$
\subsubsection{Continuation of Proof of Proposition}
Let us denote the two components of $S_\infty$ by $\mathcal{N}_1$ and $\mathcal{N}_2$ such that ${\mathcal{N}_1} \in {\mathcal{M}_{g, 1}}$ and ${\mathcal{N}_2} \in {\mathcal{M}_{0, n+1}}$. Lemma 2 says that $S_\infty$ must have at least $2g-1$ many small cuspidal eigenvalues. By \cite[Th\'{e}or\'{e}me 0.2]{O-R} the number of non-zero small eigenvalues of $\mathcal{N}_1$ is at most $2g-2$. In particular, the number of small cuspidal eigenvalues of $\mathcal{N}_1$ is at most $2g-2$. Thus for some $i$, $1 \leq i \leq 2g-1$, ${\phi^i_\infty}$ is not the zero function when restricted to ${\mathcal{N}_2}$ i.e. ${\phi^i_\infty}$ is a cuspidal eigenfunction of ${\mathcal{N}_2}$. This is a contradiction because ${\mathcal{N}_2}$ does not have any small cuspidal eigenfunction by \cite{H} or \cite{O}. $\Box$
\begin{rem}
The arguments in the proof of Proposition are applicable to more general settings. In particular, let $({S_m})$ be a sequence in $\mathcal{M}_{g, n}$ that converges to ${S_\infty} \in \partial{\mathcal{M}_{g, n}}.$ Let $({\lambda_m}, {\phi_m})$ be a normalized small eigenpair of $S_m$. Let ${\lambda_m}\rightarrow {\lambda_\infty}$ as $m$ tends to infinity. The arguments show the following: If ${\liminf_{m \rightarrow \infty}} \lVert {\phi_m} \rVert < 1$ then there exists a curve $\alpha_m$, homotopic to a geodesic of length tending to zero, on which, up to extracting a subsequence, $\phi_m$ has constant sign.
\end{rem}


\begin{thebibliography}{}
\bibitem[Bu]{Bu}   Buser, Peter, Geometry and spectra of compact Riemann surfaces.
Progress in Mathematics, 106. Birkh\"{a}user Boston, Inc., Boston, MA, 1992.

\bibitem[B-J-S]{B-J-S} Bers, L.; John, F.; Schecter, M.; Partial Differential Equations. New York: Inter-science 1964.

\bibitem[Ch]{Ch}  Cheng, Y., Eigenfunctions and nodal domains. Comment. Math. Helv. 51 (1976), 43-55

\bibitem[F]{F} Fay, J. D.; Fourier coefficient of the resolvent for a Fuchsian group. J. Reine Angew. Math. 293, 143-203 (1977).

\bibitem[H]{H}  Hejhal, D. ; Regular b-groups, degenerating Riemann surfaces and spectral theory, Memoires of Amer. Math. Soc. 88, No. 437, 1990.

\bibitem[Hu]{Hu} Huxley, M. N.; Cheeger's inequality with a boundary term, Commentarii Mathematici Helvetici 58 (1983).

\bibitem[I]{I} Iwaniec, H., Introduction to the Spectral Theory of Automorphic Forms, Bibl. Rev. Mat. Iberoamericana, Revista Matem\'{a}tica Iberoamericana, Madrid, 1995.

\bibitem[J]{J} Judge, C.; Tracking eigenvalues to the frontier of moduli space. I; J. Funct. Anal. 184 (2001), no. 2, 273- 290.

\bibitem[Ji]{Ji} Ji, Lizhen; Spectral degeneration of hyperbolic Riemann surfaces. J. Differential Geom. 38 (1993), no. 2, 263 - 313.

\bibitem[Ke]{Ke} Keen, Linda; Collars on Riemann surfaces. Discontinuous groups and Riemann surfaces, Ann. of Math. stud. 79. Princeton Univ. Press, Princeton, NJ (1974), 263-268.

\bibitem[Le]{Le} Lebedev, N. N., Special Functions and their Applications. Dover Publications, New York, 1972.

\bibitem[O]{O}   Otal, Jean-Pierre; Three topological properties of small eigenfunctions on hyperbolic surfaces. Geometry and Dynamics of Groups and Spaces, Progr. Math. 265, Birkh\"{a}user, Bassel, 2008.

\bibitem[O-R]{O-R}   Otal, Jean-Pierre; Rosas, Eulalio; Pour toute surface hyperbolique de genre g, ${\lambda_{2g-2}}>1/4$. Duke Math. J. 150 (2009), no. 1, 101 - 115.

\bibitem[P-S]{P-S} Philip, R. S. , Sarnak, P.; On cusps forms for co-finite subgroups of $PSL(2, \mathbb{R})$, Invent. Math. 80 (1985), 339 - 364.

\bibitem[W]{W} Wolf, M.; Infinite energy harmonic maps and degeneration of hyperbolic surfaces in Moduli spaces, J. Differential Geometry 33 (1991) 487-539

\bibitem[Wo]{Wo} Wolpert, S. A.; Spectral limits for hyperbolic surface, I, Invent. Math. 108 (1992), 67 - 89
\end{thebibliography}
\end{document}